\documentclass{amsart}
\numberwithin{equation}{section}
\usepackage[english]{babel}   \usepackage[T1]{fontenc}
\usepackage[latin1]{inputenc}
\usepackage{indentfirst}
\usepackage{enumitem}
\usepackage{amsmath,amssymb, amsbsy}
\usepackage{amsfonts}
\usepackage{hyperref}
\usepackage{latexsym}
\usepackage{amsthm}
\usepackage[dvips]{graphicx}
\DeclareGraphicsExtensions{.pdf,.png,.jpg,.eps}

\usepackage{verbatim}

\usepackage[usenames,dvipsnames]{xcolor}
\usepackage{color}
\usepackage[latin1]{inputenc}
\usepackage[active]{srcltx}
\definecolor{Arancio}{cmyk}{0,0.61,0.87,0}

\newcommand{\brd}[1]{\mathbb{#1}}
\newcommand{\R}{\brd{R}}
\newcommand{\abs}[1]{\left\lvert {#1} \right\rvert}
\newcommand{\norm}[2]{\left\Vert {#1} \right\Vert_{#2}}

\newtheorem{teo}{Theorem}[section]
\newtheorem{Corollary}[teo]{Corollary}
\newtheorem{Lemma}[teo]{Lemma}
\newtheorem{Theorem}[teo]{Theorem}
\newtheorem{Proposition}[teo]{Proposition}
\theoremstyle{definition}

\newtheorem{remark}[teo]{Remark}
\newtheorem{Assumption}[teo]{Assumption}
\newtheorem{Notation}[teo]{Notation}

\begin{document}

\title[Regularity for odd solutions to degenerate or singular problems]
{Liouville type theorems and regularity of solutions to degenerate or singular problems part II: odd solutions}

\author{Yannick Sire, Susanna Terracini and Stefano Vita}

\address[Y. Sire]{Department of Mathematics
	\newline\indent
	Johns Hopkins University
	\newline\indent
	3400 N. Charles Street, Baltimore, MD 21218, U.S.A.}
\email{sire@math.jhu.edu}

\address[S. Terracini]{Dipartimento di Matematica G. Peano
	\newline\indent
	Universit\`a degli Studi di Torino
	\newline\indent
	 Via Carlo Alberto 10, 20123 Torino, Italy}
\email{susanna.terracini@unito.it}

\address[S. Vita]{Dipartimento di Matematica
	\newline\indent
	Universit\`a degli Studi di Milano Bicocca
	\newline\indent
	 Piazza dell'Ateneo Nuovo 1, 20126, Milano, Italy}
\email{stefano.vita@unimib.it}

\date{\today}
\thanks{Work partially supported by the ERC Advanced Grant 2013 n.~339958 Complex Patterns for Strongly Interacting Dynamical Systems - COMPAT. Y.S. was partially supported by the Simons foundation.}



\maketitle

\begin{center}{\it To Sandro, with friendship, admiration and much more}\end{center}

\begin{abstract}
We consider a class of equations in divergence form with a singular/degenerate weight \[
-\mathrm{div}(|y|^a A(x,y)\nabla u)=|y|^a f(x,y)+\textrm{div}(|y|^aF(x,y))\;.
\]
Under suitable regularity assumptions for the matrix $A$, the forcing term $f$ and the field $F$, we prove  H\"older continuity of solutions which are odd in $y\in\mathbb{R}$, and possibly of their derivatives. In addition, we show stability of the $C^{0,\alpha}$ and $C^{1,\alpha}$ a priori bounds for approximating problems in the form
\[
-\mathrm{div}((\varepsilon^2+y^2)^{a/2} A(x,y)\nabla u)=(\varepsilon^2+y^2)^{a/2} f(x,y)+\textrm{div}((\varepsilon^2+y^2)^{a/2}F(x,y))
\]
as $\varepsilon\to 0$. Our method is based upon blow-up and appropriate Liouville type theorems.

\end{abstract}

\tableofcontents

\section{Introduction and main results}
Let $z=(x,y)\in\mathbb{R}^{n+1}$, with $x\in\mathbb{R}^{n}$ and $y\in\mathbb{R}$, $n\geq1$, $a\in\mathbb{R}$. Our aim is to study the boundary behaviour of solutions to a class of problems involving singular/degenerate operators in divergence form including
$$\mathcal L_au:=\mathrm{div}(|y|^aA(x,y)\nabla u),$$
and their regularizations. The boundary here coincides with $\Sigma:=\{y=0\}$ the \textsl{characteristic manifold}, where the weight becomes degenerate or singular, and this happens respectively when $a>0$ and $a<0$. Accordingly, this class of operators is called degenerate elliptic.

The first motivation for this work is to complete the study started in \cite{SirTerVit1} on local regularity for solutions to degenerate/singular problems including the following
\begin{equation}\label{La1}
-\mathrm{div}\left(|y|^a\nabla u\right)=|y|^af+\mathrm{div}\left(|y|^aF\right)\qquad\mathrm{in \ } B_1.
\end{equation}
In \cite{SirTerVit1}, we treated the regularity of even-in-$y$ solutions (corresponding to Neumann boundary conditions), including the case of variable coefficients. We provided local $\mathcal C^{0,\alpha}$  and $\mathcal C^{1,\alpha}$ estimates, which are uniform as the parameter $\varepsilon\to0^+$, for  even solutions of regularized uniformly elliptic problems of the form
\begin{equation}\label{LarhoA}
-\mathrm{div}\left(\rho_\varepsilon^a(y)A(x,y)\nabla u_\varepsilon\right)=\rho_\varepsilon^a(y)f_\varepsilon+\mathrm{div}\left(\rho_\varepsilon^a(y)F_\varepsilon\right)\qquad\mathrm{in \ } B_1,
\end{equation}
 where the regularized family of weights $\rho_\varepsilon^a$ is defined as: 
\begin{equation}\label{rho}
\rho_\varepsilon^a(y):=\begin{cases}
(\varepsilon^2+y^2)^{a/2}\min\{\varepsilon^{-a},1\} &\mathrm{if \ }a\geq0,\\
(\varepsilon^2+y^2)^{a/2}\max\{\varepsilon^{-a},1\} &\mathrm{if \ }a\leq0.
\end{cases}
\end{equation}

A further motivation comes from a remarkable link between our operators and fractional powers of the Laplacian, from a Dirichlet-to-Neumann point of view, as highlighted in \cite{CafSil1},  when our weights belong to the $A_2$-class; i.e. $a\in(-1,1)$.

Goal of this paper is to deal with odd-in-$y$  solutions to \eqref{La1} (corresponding to Dirichlet boundary conditions), providing local regularity, when possible in the $\varepsilon$-stable sense, by proving uniform bounds for solutions to \eqref{LarhoA}.
Odd solutions make sense as energy solutions in the natural weighted Sobolev spaces whenever $a\in(-\infty,1)$ (in the sense of \S\ref{sect:sobolev}). At first, we notice that can not expect, for the odd solutions,  the same estimates as for the even ones, where the regularity results from the combined effect of the ellipticity and the boundary condition.  In fact, the function $y|y|^{-a}$ is $\mathcal L_a$-harmonic with finite energy when $a<1$ (in case of $A=\mathbb I$), and for $a\in(0,1)$ is no more than H\"older continuous. We will refer to this special solution as the \textsl{characteristic odd comparison solution}. Similar, yet smoother, characteristic odd comparison solutions exist for the full regularized family of $\varepsilon$-problems (in a rather general setting). Nonetheless, one major obstruction in the study of regularity is the fact that our weights need not to be locally integrable when $a\leq-1$, preventing the application of classical regularity theory such as that developed for degenerate weights of the $A_2$-Muckenhoupt class, starting from the seminal papers \cite{FabKenSer,FabJerKen1,FabJerKen2}. We shall adopt here a different perspective, exploiting suitably tailored Liouville type theorems as main tools.  To this aim, a major hindrance is that the measure $|y|^a\mathrm{d}z$ is not absolutely continuous with respect to the Lebesgue measure. In order to overcome this difficulty, one can be guided by the following insight:
\begin{Proposition}\label{prop1}
Let $a\in(-\infty,1)$ and $u\in H^{1,a}(B_1)$ be an odd energy solution to \eqref{La1} in $B_1$ (for simplicity with $F=0$). Then for any $r<1$ the ratio $w=u/y|y|^{-a}\in H^{1,2-a}(B_r)$ and it is an even energy solution to
\begin{equation}\label{BHLaw}
-\mathrm{div}\left(|y|^{2-a}\nabla w\right)=|y|^{2-a}\bar f=|y|^{2-a}\frac{f}{y|y|^{-a}}\qquad\mathrm{in \ }B_r.
\end{equation}
\end{Proposition}
Proposition \ref{prop1} allows the application of the results for even solutions already proved in \cite{SirTerVit1},  providing regularity up to the multiplicative factor $y|y|^{-a}$. Thanks to this observation it is natural to shift the study of regularity for odd solutions to that of even solutions of  the auxiliary problem above. A similar perspective has been adopted in \cite{ShaYer} for the obstacle problem in the same singular/degenerate setting.

As an example,  by the Schauder estimates in \cite{SirTerVit1},  when the forcing $\bar f=:f/y|y|^{-a}$  in \eqref{BHLaw} is $C^{k,\alpha}$, then the ratio $w=u/y|y|^{-a}$ is locally $C^{k+2,\alpha}$.  Thus, we understand that the correct way to face the regularity of odd solutions consists in seeking $\mathcal C^{0,\alpha}$  and $\mathcal C^{1,\alpha}$ bounds for the ratio between the solution and the characteristic odd one, depending on the regularity of the same ratio of the right hand side. This point of view corresponds to \textsl{(possibly higher order and/or non homogeneous) boundary Harnack principle} at $\Sigma$ in the sense of \cite{CafFabMorSal, DesSav, FabJerKen2, JavNeu, JerKen}. It is worthwhile noticing  that, when $a\in(-\infty,1)$, then the exponent $2-a$ belongs to $(1,+\infty)$, placing equation \eqref{BHLaw} in the so called super degenerate case, again outside the land of $A_2$-Muckenhoupt weights theory, and which has been treated in \cite{SirTerVit1} when associated with Neumann boundary conditions. Furthermore, looking at the right hand side of \eqref{BHLaw}, we realize that the transition from the odd to the even case requires to pay a cost in terms of more stringent conditions on the forcing term $f$, in the sense that the ratio $\frac{f}{y|y|^{-a}}$ must possess some regularity (integrability at least); in other words, when $a<0$, it means that the forcing term is vanishing with a certain rate at $\Sigma$. In this regard, our results are connected with the recent paper \cite{AllLSha}, where a boundary Harnack principle with right hand side is established in the uniformly elliptic case.

As already pointed out, our results are not limited to the $A_2$-Muckenhoupt class of weights, which restricts $a$ in the interval $(-1,1)$. Nonetheless, we wish to state the following corollary, which joins the results contained in this paper with the Schauder theory for even solutions developed in \cite{SirTerVit1}, concerning full regularity for energy solutions of degenerate or singular problems when the weight is $A_2$-Muckenhoupt and $A=\mathbb I$.
\begin{Corollary}
Let $a\in(-1,1)$, $k\in\mathbb N\cup\{0\}$, $\alpha\in(0,1)$ and consider $u\in H^{1,a}(B_1)$ an energy solution to
\begin{equation*}
-\mathrm{div}\left(|y|^a\nabla u\right)=|y|^af\qquad\mathrm{in \ }B_1.
\end{equation*}
Let us consider the even and odd parts\footnote{Even and odd parts (in $y$) of a function are defined as usual as
$$f_e(x,y)=\frac{f(x,y)+f(x,-y)}{2},\qquad f_o(x,y)=\frac{f(x,y)-f(x,-y)}{2}.$$}
 (with respect to $y$) of the forcing term $f$. Let
$$f=f_e+f_o=f_e+y|y|^{-a}\tilde f_e\qquad\mathrm{with \ }f_e,\tilde f_e\in C^{k,\alpha}(B_1).$$
Then
\begin{equation*}
u=u_e+u_o=u_e+y|y|^{-a}\tilde u_e,\qquad\mathrm{with \ }u_e,\tilde u_e\in C^{k+2,\alpha}_{\mathrm{loc}}(B_1).
\end{equation*}
\end{Corollary}

As a next step, we aim at deepening the $\varepsilon$-stability of these estimates with respect to the family of regularized weights \eqref{rho} (also including the variable coefficient case). In other words, we deal with odd-in-$y$ solutions to the family of equations \eqref{LarhoA}.
We will provide local uniform-in-$\varepsilon$ regularity estimates, enlightening their delicate link with curvature issues related with the matrix $A$. As we shall see, also the notion of characteristic solution must be suitably  adjusted in order to deal with the variable coefficient cases.  Finally, we will apply our results to a family of degenerate/singular equations naturally associated with the euclidean Laplacian expressed in Fermi coordinates in the neighbourhood of an embedded hypersurface. 

Below we set the minimal assumptions on the matrix $A$ that we need throughout the paper:
\begin{Assumption}[HA]\label{(HA)}
The matrix $A=(a_{ij})$ is $(n+1,n+1)$-dimensional and symmetric $A=A^T$, has the following symmetry with respect to $\Sigma$: 
we have
\begin{equation*}
A(x,y)=JA(x,-y)J,\qquad\qquad\mathrm{with}\qquad\qquad
J
=\left(
\begin{array}{c|c}
\mathbb I_n & 0 \\
\hline
0 & -1
\end{array}
\right).
\end{equation*}
Therefore, $A$ is continuous and satisfies the uniform ellipticity condition $\lambda_1|\xi |^2 \leq A(x,y)\xi\cdot\xi \leq \lambda_2|\xi |^2$, for all $\xi\in\R^{n+1}$, for every $(x,y)$ and some ellipticity constants $0<\lambda_1\leq\lambda_2$. Moreover, the characteristic manifold $\Sigma$ is assumed to be invariant with respect to $A$; that is, there exists a suitable scalar function $\mu$ such that there exists a positive constant such that
\begin{equation}\label{boundmu}
\frac{1}{C}\leq\mu(x,y)\leq C
\end{equation}
and with
$$A(x,0)\cdot e_{y}=\mu(x,0) e_y.$$
\end{Assumption}
Whenever the hypothesis on $A$ are not specified, we always imply Assumption (HA). From now on, through out the paper, whenever not otherwise specified, in order to simplify the notations, we will work with $A=\mathbb I$ every time this condition is not playing a role in the proofs. In the perspective of Proposition \ref{prop1}, but considering odd solutions for the family of regularized problems in \eqref{LarhoA}, it will be convenient to adopt the following notation on the matrix $A$.
\begin{Notation}[HA+]\label{(HA+)}
We can always write matrix $A$ as:
\begin{equation*}
A(x,y)=\mu(x,y)B(x,y),
\end{equation*}
with
\begin{equation*}
B(x,y)=
\left(
\begin{array}{c|c}
\tilde B(x,y) & T(x,y) \\
\hline
T(x,y) & 1
\end{array}
\right),
\end{equation*}
where $\tilde B$ is a $(n,n)$-dimensional matrix and $T:\R^{n+1}\to\R^n$ (we denote by $\tilde A=\mu\tilde B$). We remark here that under our hypothesis on the symmetries of coefficients; one has, for $y<0$
\begin{equation*}
A(x,y)=\mu(x,-y)\left(
\begin{array}{c|c}
\tilde B(x,-y) & -T(x,-y) \\
\hline
-T(x,-y) & 1
\end{array}
\right).
\end{equation*}
\end{Notation}

The structural assumption on the matrix $A$ is consistent with \cite{CafSti}. Moreover, it fits also with the metric induced by Fermi's coordinates, which allow to study phenomena of singularity or degeneration on a characteristic manifold $\Sigma$ which is a generic (regular enough) $n$-dimensional hypersurface embedded in $\R^{n+1}$.
Hence, the objective will be to consider the ratio $w_\varepsilon$ between odd solutions $u_\varepsilon$ to \eqref{LarhoA} and functions of the form
\begin{equation}\label{veps}
v^a_\varepsilon(x,y)=(1-a)\int_0^y\rho_\varepsilon^{-a}(s)\mu(x,s)^{-1}\mathrm{d}s,
\end{equation}
which play now the role of the \textsl{characteristic odd solution} for the regularized family of weights in the variable coefficients case. It is worthwhile stressing that the characteristic solutions $v^a_\varepsilon$ do not longer solve the homogenous problem, as a dependence on the curvature appears. 

As said, we wish to obtain uniform local regularity estimates for $w_\varepsilon$ which will be even solutions to an auxiliary weighted problems having the following structure
\begin{equation}\label{simpleBH}
-\mathrm{div}\left(\rho_\varepsilon^a(v_\varepsilon^a)^2A\nabla w_\varepsilon\right)=\rho_\varepsilon^a(v_\varepsilon^a)^2f_\varepsilon+\mathrm{div}\left(\rho_\varepsilon^a(v_\varepsilon^a)^2F_\varepsilon\right)+\rho_\varepsilon^a(v_\varepsilon^a)^2b_\varepsilon\cdot\nabla w_\varepsilon.
\end{equation}
The new weights appearing in the auxiliary equation are equivalent, though not equal, (using \eqref{boundmu}) to
\begin{equation}
\omega_\varepsilon^a(y)
=\rho_\varepsilon^a(y)(1-a)^2(\chi_\varepsilon^a(y))^2
\end{equation}
where we have defined
\begin{equation}\label{eq:chi}
\chi_\varepsilon^a(y):=\int_0^y\rho_\varepsilon^{-a}(s)\mathrm{d}s\;.
\end{equation}
We remark that, as $a\in(-\infty,1)$, such a class of weights is always super degenerate; indeed, at $\Sigma$, they behave like
$$\omega_\varepsilon^a(y)\sim\begin{cases}
y^2 &\mathrm{if \ }\varepsilon>0\\
|y|^{2-a} &\mathrm{if \ }\varepsilon=0,
\end{cases}$$
with $2-a\in(1,+\infty)$.

Our first main result concerns in fact the even solutions to the auxiliary family of equations \eqref{simpleBH}. It essentially consists in extending (in a non trivial way) the analogous result already obtained in \cite{SirTerVit1} to the new family of weights $\rho_\varepsilon^a(v_\varepsilon^a)^2$.
\begin{Theorem}\label{theo:first}
Let $a\in(-\infty,1)$ and, as $\varepsilon\to0$, let $\{w_\varepsilon\}$ be a family of solutions in $B_1^+$ of \eqref{simpleBH} which are even-in-$y$; that is, satisfying the boundary condition
$$\rho_\varepsilon^a(v_\varepsilon^a)^2\partial_yw_\varepsilon=0 \qquad \mathrm{on \ }\partial^0B_1^+.$$
$1)$ Let $r\in(0,1)$, $\beta>1$, $p_1>\frac{n+3+(-a)^+}{2}$, $p_2,p_3>n+3+(-a)^+$, and $\alpha\in(0,2-\frac{n+3+(-a)^+}{p_1}]\cap(0,1-\frac{n+3+(-a)^+}{p_2}]\cap(0,1-\frac{n+3+(-a)^+}{p_3}]$. Let's moreover take $A$ with continuous coefficients and $\| b_\varepsilon\|_{L^{p_3}(B_1^+,\omega_\varepsilon^a(y)\mathrm{d}z)}\leq b$. There is a positive constant $c$ depending on $a$, $b$, $n$, $\beta$, $p_1$, $p_2$, $p_3$, $\alpha$ and $r$ only such that functions $w_\varepsilon$
satisfy
\begin{equation*}
\|w_\varepsilon\|_{C^{0,\alpha}(B_r^+)}\leq c\left(\|w_\varepsilon\|_{L^\beta(B_1^+,\omega_\varepsilon^a(y)\mathrm{d}z)}+ \| f_\varepsilon\|_{L^{p_1}(B_1^+,\omega_\varepsilon^a(y)\mathrm{d}z)}+\|F_\varepsilon\|_{L^{p_2}(B_1^+,\omega_\varepsilon^a(y)\mathrm{d}z)}\right).
\end{equation*}
$2)$ Let $r\in(0,1)$, $\beta>1$, $p_1,p_2>n+3+(-a)^+$, and $\alpha\in(0,1-\frac{n+3+(-a)^+}{p_1}]\cap(0,1-\frac{n+3+(-a)^+}{p_2}]$. Let $F_\varepsilon=(F^1_\varepsilon,...,F^{n+1}_\varepsilon)$ with the $y$-component vanishing on $\Sigma$: $F^{n+1}_\varepsilon(x,0)= F^y_\varepsilon(x,0)=0$ in $\partial^0B_1^+$.
Let's moreover take $A$ with $\alpha$-H\"older continuous coefficients and $\|b_\varepsilon\|_{L^{2p_2}(B_1^+,\omega_\varepsilon^a(y)\mathrm{d}z)}\leq b$. There is a positive constant $c$ depending on $a$, $b$, $n$, $\beta$, $p_1$, $p_2$, $\alpha$ and $r$ only such that functions $w_\varepsilon$ satisfy
\begin{equation*}
\|w_\varepsilon\|_{C^{1,\alpha}(B_r^+)}\leq c\left(\|w_\varepsilon\|_{L^\beta(B_1^+,\omega_\varepsilon^a(y)\mathrm{d}z)}+ \|f_\varepsilon\|_{L^{p_1}(B_1^+,\omega_\varepsilon^a(y)\mathrm{d}z)}+\|F_\varepsilon\|_{C^{0,\alpha}(B_1^+)}\right).
\end{equation*}
\end{Theorem}
We would like to remark here that local $C^{2,\alpha}$ uniform estimates (up to $\Sigma$) with respect to the regularization can not be proven (for a counterexample we refer to \cite[Remark 5.4]{SirTerVit1}).

When applying Theorem \ref{theo:first} to the quotient
\begin{equation}
w_\varepsilon=\dfrac{u}{v_\varepsilon^a}
\end{equation}
of a solution of \eqref{LarhoA} and the characteristic solution \eqref{veps}, we realise that the actual terms appearing in right hand side of \eqref{simpleBH} depend on the original forcings $f,F$ jointly with the parameters $\mu,T,B$ of the matrix  $A$ written as in Notation (HA+). In particular, as shown in \eqref{BHepseq1}, we see the appearance of a drift term involving the $x$-derivatives of $\mu$ which, consequently, need to satisfy a $\mathcal C^{0,\alpha}$ condition.  Our main result is Theorem \ref{holderBH}. 
We give here below a simplified statement, suitable to be applied to the case of  laplacians in Fermi coordinates treated in subsection \S\ref{subsect:fermi}.

\begin{Theorem}\label{holderBHsimple}
Let $a\in(-\infty,1)$, the matrix $A$ written as in Notation (HA+) with $T\equiv 0$. As $\varepsilon\to0$ let $\{u_\varepsilon\}$ be a family of solutions in $B_1^+$ of
\begin{equation*}\label{1oddBH}
\begin{cases}
-\mathrm{div}\left(\rho_\varepsilon^aA\nabla u_\varepsilon\right)=\rho_\varepsilon^af_\varepsilon+\mathrm{div}\left(\rho_\varepsilon^aF_\varepsilon\right) & \mathrm{in \ } B_1^+\\
u_\varepsilon=0 & \mathrm{on \ }\partial^0B_1^+.
\end{cases}
\end{equation*}
Let also $\{v_\varepsilon^a\}$ be the family of solutions defined in \eqref{veps} in $B_1^+$.  Denote
$$w_\varepsilon=\frac{u_\varepsilon}{v_\varepsilon^a}\;.$$

$1)$ Assume $\mu$ be Lipschitz continuous, $r\in(0,1)$, $\beta>1$, $p_1>\frac{n+3+(-a)^+}{2}$, $p_2>n+3+(-a)^+$, and $\alpha\in(0,2-\frac{n+3+(-a)^+}{p_1}]\cap(0,1-\frac{n+3+(-a)^+}{p_2}]$. Let's moreover take $A$ with continuous coefficients. There is a positive constant $c$ depending on $a$, $n$, $\beta$, $p_1$, $p_2$, $\alpha$ and $r$ only such that the $w_\varepsilon$ satisfy
\begin{eqnarray*}
\|w_\varepsilon\|_{C^{0,\alpha}(B_r^+)} &\leq& c\left( \|w_\varepsilon\|_{L^\beta(B_1^+,\omega_\varepsilon^a(y)\mathrm{d}z)}+ \|{f_\varepsilon}/{v_\varepsilon^a}\|_{L^{p_1}(B_1^+,\omega_\varepsilon^a(y)\mathrm{d}z)} \right.\\
&& \left. + \|F_\varepsilon^y/(yv_\varepsilon^a)\|_{L^{p_1}(B_1^+,\omega_\varepsilon^a(y)\mathrm{d}z)} +\|F_\varepsilon/v_\varepsilon^a\|_{L^{p_2}(B_1^+,\omega_\varepsilon^a(y)\mathrm{d}z)}\right).
\end{eqnarray*}
$2)$ Assume $\mu\in \mathcal C^{1,\alpha}(B_1^+)$, and let $r\in(0,1)$, $\beta>1$, $p_1>n+3+(-a)^+$, and $\alpha\in(0,1-\frac{n+3+(-a)^+}{p_1}]$. Let $F_\varepsilon=(F^1_\varepsilon,...,F^{n+1}_\varepsilon)$ with the $\alpha$-H\"older continuous ratio between the $y$-component and $v_\varepsilon^a$ vanishing on $\Sigma$: $F^{n+1}_\varepsilon(x,0)/v_\varepsilon^a= F^y_\varepsilon(x,0)/v_\varepsilon^a=0$ in $\partial^0B_1^+$. Let's moreover take $A$ with $\alpha$-H\"older continuous coefficients. There is a positive constant $c$ depending on $a$, $n$, $\beta$, $p_1$, $\alpha$ and $r$ only such that 
\begin{eqnarray*}
\|w_\varepsilon\|_{C^{1,\alpha}(B_r^+)} &\leq& c\left( \|w_\varepsilon\|_{L^\beta(B_1^+,\omega_\varepsilon^a(y)\mathrm{d}z)}+ \|{f_\varepsilon}/{v_\varepsilon^a}\|_{L^{p_1}(B_1^+,\omega_\varepsilon^a(y)\mathrm{d}z)} \right.\\
&& \left. + \|F_\varepsilon^y/(yv_\varepsilon^a)\|_{L^{p_1}(B_1^+,\omega_\varepsilon^a(y)\mathrm{d}z)} +\|F_\varepsilon/v_\varepsilon^a\|_{C^{0,\alpha}(B_1^+)}\right).
\end{eqnarray*}
\end{Theorem}
It is worthwhile  noticing here that any energy odd solution to \eqref{LarhoA} for $\varepsilon=0$ (under suitable conditions on the matrix and the right hand side) can be approximated by a $\varepsilon$-sequence of solutions to \eqref{LarhoA} satisfying the hypothesis in our regularity results. The same happens for the auxiliary weighed problems solved by the even functions $w=u/y|y|^{-a}$. This is done in details in \cite[Section 2 and 6]{SirTerVit1}.
\begin{remark}
A special, yet fundamental, case is when take $A=\mathbb I$, so that $\mu\equiv 1$ and the family of fundamental comparison odd solutions $v_\varepsilon^a$'s are in fact the $\chi_\varepsilon^a$'s. Nevertheless, it has to be noticed that, in the presence of non trivial curvature, the ratio $v_\varepsilon^a/\chi_\varepsilon^a$ may not be uniformly (in $\varepsilon$) bounded in $\mathcal C^{1,\alpha}(B_1^+)$. Furthermore, in the variable coefficient case, the $\chi_\varepsilon^a$'s are not in the kernel of the corresponding operators, as a (possibly weird) right hand side appears.
\end{remark}

This Theorem  finds a natural application to the study of  the boundary behaviour of solutions of operators degenerate/singular at embedded manifolds, as shown by the following result.

\begin{Corollary}\label{holderBHFermi}
Let $\Sigma$ be an $n$-dimensional hypersurface embedded in $\R^{n+1}$, of class $\mathcal C^{3,\alpha}$ and let $d_\Sigma(X)$ denote the signed distance of $X$ to $\Sigma$. Let $a\in(-\infty,1)$, $R>0$ sufficiently small, and consider, as $\varepsilon\to0$,  a family of solutions to 
\begin{equation*}\label{1oddBHsigma}
\begin{cases}
-\mathrm{div}\left(\rho_\varepsilon^a\circ d_\Sigma\nabla u_\varepsilon\right)=\rho_\varepsilon^a\circ d_\Sigma f_\varepsilon+\mathrm{div}\left(\rho_\varepsilon^a\circ d_\Sigma F_\varepsilon\right) & \mathrm{in \ } B_R\cap\{d_\Sigma(X)>0\}\\
u_\varepsilon=0 & \mathrm{on \ } B_R\cap \Sigma.
\end{cases}
\end{equation*}
Let also $\{\chi_\varepsilon^a\}$ be the family of functions defined in \eqref{eq:chi} in $B_R$.  Denote
$$w_\varepsilon=\frac{u_\varepsilon}{\chi_\varepsilon^a\circ d_\Sigma }\;,$$
$1)$ The same conclusion of point 1) of Theorem \ref{holderBHsimple} holds with $v_\varepsilon^a$ replaced by $\chi_\varepsilon^a$, $y$ by $d_\Sigma(X)$ and $e_{n+1}$  by the normal $\nu$ at $\Sigma$.\\
$2)$ The same conclusion of point 2) of Theorem \ref{holderBHsimple} holds in $C^{1,\alpha}(B_r\cap\{y\geq \sqrt{\varepsilon}\})$ where $c$ is independent of $\varepsilon$, and,  again, $v_\varepsilon^a$ replaced by $\chi_\varepsilon^a$, $y$ by $d_\Sigma(X)$ and $e_{n+1}$  by the normal $\nu$ at $\Sigma$.\\
\end{Corollary}

\begin{remark}
In particular, letting $\varepsilon\to 0$ we find $C^{1,\alpha}(B_r^+)$ estimates in the degenerate/singular case, though not in the full $\varepsilon$-stable sense. The reason is the possible lack of uniform-in-$\varepsilon$ smoothness of the ratio $v_\varepsilon^a/\chi_\varepsilon^a$.
\end{remark}

\subsection{Proof of Corollary \ref{holderBHFermi}}\label{subsect:fermi}
As already mentioned, the structural assumption on the matrix $A$ done in Assumption (HA) with Notation (HA+) fits also with the metric induced by Fermi's coordinates around the characteristic manifold $\Sigma$ (see \cite{PacWei}). Let $\Sigma$ be an oriented regular enough hypersurface embedded in $\R^{n+1}$. We are concerned with operators associated with Dirichlet energies of the form
\begin{equation*}
\int_{\{d_\Sigma(X)>0\}}(\rho_\varepsilon^a\circ d_\Sigma)(X)|\nabla u|^2,
\end{equation*}
with $a\in\R$, 
$X \in \R^{n+1}$ and $d_\Sigma(\cdot)$ the signed distance function to $\Sigma$. Let $g_e$ be the Euclidean metric on $\R^{n+1}$ and denote by $\nu$ the unit normal vector field on $\Sigma$.  We define Fermi coordinates in a tubular neighborhood of  $\Sigma$ as follows: let $z\in\Sigma$ and $y\in\R$, and define
\begin{equation*}
Z(z,y):=z+y\,\nu(z).
\end{equation*}
Nevertheless, $z$ and $y\,\nu(z)$ belong to $\R^{n+1}$. Given $y \in \R$, we define  
$$\Sigma_y:=\{Z(z,y) \in \R^{n+1}\ : \ z\in\Sigma\}. $$
Following Lemma 6.1 in \cite{PacWei}, one has that the induced metric on $\Sigma_y$ is given by
\begin{equation*}
g^y= g^0-2yh^0+y^2h^0\otimes h^0,
\end{equation*}
where $g^0$ is the induced metric on $\Sigma$, $h^0$ is the second fundamental form on $\Sigma$ and $h^0\otimes h^0$ its square; namely we have 
$$
h^0(t_1,t_2)=-g^0(\nabla_{t_1}\nu,t_2)
$$
for all $t_1,t_2$ on the tangent bundle of $\Sigma$. Notice in particular that, in local coordinates, the terms $g^0,h^0,h^0\otimes h^0$ depend only on $z$. 
Therefore, invoking Lemma 6.3 in \cite{PacWei}, one finally has 
$$
Z^*g_e=g^y+dy^2
$$
where $g^y$ is considered as a family of metrics on the tangent bundle of $\Sigma$, depending smoothly on $y$ in a neighborhood of $0$ in $\R$.

In other words, we are obtaining a quadratic form for $v(z,y)=u(Z(z,y))$ of the form
\begin{equation*}
\int_{0}^{y_0} \rho_\varepsilon^a(y)\int_{\Sigma_y}\left(|\nabla_{g^y}v|^2+|\partial_yv|^2\right) \sqrt{\mathrm{det} g^y}.
\end{equation*}

Recall that the variation with respect to $y$ of of the volume form of the parallel hypersurfaces $\Sigma_y$ satisfy the equation:
\begin{equation}\label{eq:gy}
H_y=-\dfrac{1}{\sqrt{\mathrm{det} g^y}}\dfrac{d}{dy}\sqrt{\mathrm{det} g^y}.
\end{equation}
Hence, by considering a parametrization of $\Sigma$ of the form $z=\psi(x)$ with $x\in \R^n$, then one obtains for $w(x,y)=v(\psi(x),y)$
\begin{equation*}
\int \rho_\varepsilon^a(y)A\nabla w\cdot\nabla w,
\end{equation*}
where
\begin{equation*}
A(x,y)=
\left(
\begin{array}{c|c}
\tilde A (x,y) & 0 \\
\hline
0 & 1
\end{array}
\right)\cdot\sqrt{\det g^y}.
\end{equation*}
We remark that the matrix $A$ satisfies Assumption (HA), and can be expressed as in Notation (HA+) with $\mu(x,y)=\sqrt{\det g^y}$. As $\Sigma\in\mathcal C^{3,\alpha}$, we have $\mu\in\mathcal C^{1,\alpha}(B_{r_0}^+)$ for $r_0$ small enough. Hence we are in the position to apply Theorem \ref{holderBHsimple}. Next we have to compare the two families $v_\varepsilon^a$ and $\chi_\varepsilon^a$.  At first, in order to prove point 1) we remark that Proposition \ref{c0alphaG} ensures uniform-in-$\varepsilon$ $\mathcal C^{0,\alpha}$ estimates for the ratio $v_\varepsilon^a/\chi_\varepsilon^a$. 
Using \eqref{eq:gy}, we infer that $\partial_y\mu(\cdot,0)\in\mathcal C^{1,\alpha}(B_{r_0}^+)$ and finally, by virtue of Proposition \ref{c1alphaG}, we obtain that also the ratio $v_\varepsilon^a/\chi_\varepsilon^a$ satisfies the desired uniform bounds in  $\mathcal C^{1,\alpha}(B_r\cap\{y\geq\sqrt{\varepsilon}\})$, for $r<r_0$.

\subsection*{\sl Notations.}
Below is the list of symbols we shall use throughout this paper.\\
 
\begin{tabular}{ll}
$\R^{n+1}_+=\R^n\times(0,+\infty)$ & $z=(x,y)$ with $x\in\R^n$, $y>0$ \\
$\Sigma=\{y=0\}$ & characteristic manifold \\
$B_r^+=B_r\cap\{y>0\}$ & half ball \\
$\partial^+B_r^+=S^n_+(r)=\partial B_r\cap\{y>0\}$ & upper boundary of the half ball \\
$\partial^0B_r^+=B_r\cap\{y=0\}$ & flat boundary of the half ball \\
$\rho_\varepsilon^a(y)=\left(\varepsilon^2+y^2\right)^{a/2}$ & regularized weight \\
$\omega_\varepsilon^a(y)=\rho_\varepsilon^a(y)\pi_\varepsilon^a(y)$ & regularized auxiliary weight \\
$\mathcal L_{\rho_\varepsilon^a}u=\mathrm{div}\left(\rho_\varepsilon^a(y)A(x,y)\nabla u\right)$ & regularized operator \\
$H^{1}(\Omega,\rho_\varepsilon^a(y)\mathrm{d}z)$ & weighted Sobolev space given by the completion of $C^\infty(\overline\Omega)$ \\
$H^{1}_0(\Omega,\rho_\varepsilon^a(y)\mathrm{d}z)$ & weighted Sobolev space given by the completion of $C^\infty_c(\Omega)$ \\
$\tilde H^{1}(\Omega,\rho_\varepsilon^a(y)\mathrm{d}z)$ & weighted Sobolev space given by the completion of $C^\infty_c(\overline\Omega\setminus\Sigma)$ \\
$H^{1,a}(\Omega)=H^{1}(\Omega,|y|^a\mathrm{d}z)$ & weighted Sobolev space for $\varepsilon=0$ \\
$\partial_y^au=|y|^a\partial_yu$ & "weighted" derivative \\
$y|y|^{-a}$ & characteristic odd solution \\
$v_\varepsilon^a$ & characteristic odd solution in presence of $A$ and $\varepsilon>0$\\
$a^+=\max\{a,0\}$ &\\
\end{tabular}

\section{Functional setting and preliminary results}\label{sect:sobolev}

In this section we collect the natural notions of Sobolev spaces, and their main properties, needed to work in our degenerate or singular context (for further details see \cite[Section 2]{SirTerVit1}). Let $\Omega\subset\R^{n+1}$ be non empty, open and bounded. Denoting by $C^\infty(\overline\Omega)$ the set of real functions $u$ defined on $\overline\Omega$ such that the derivatives $D^\alpha u$ can be continuously extended to $\overline\Omega$ for all multiindices $\alpha$, then for any $a\in\mathbb{R}$, $\varepsilon\geq0$ we define the weighted Sobolev space $H^1(\Omega,\rho_\varepsilon^a(y)\mathrm{d}z)$ as the closure of $C^\infty(\overline\Omega)$ with respect to the norm
$$\|u\|_{H^{1}(\Omega,\rho_\varepsilon^a(y)\mathrm{d}z)}=\left(\int_{\Omega}\rho_\varepsilon^au^2+\int_{\Omega}\rho_\varepsilon^a|\nabla u|^2\right)^{1/2}.$$
To simplify the notation we will denote 
$$H^{1,a}(\Omega)=H^{1}(\Omega,|y|^a\mathrm{d}z)=H^{1}(\Omega,\rho_0^a(y)\mathrm{d}z).$$
In the same way, we define $H^{1}_0(\Omega,\rho_\varepsilon^a(y)\mathrm{d}z)$ as the closure of $C^\infty_c(\Omega)$ with respect to the norm
$$\|u\|_{H^{1}_0(\Omega,\rho_\varepsilon^a(y)\mathrm{d}z)}=\left(\int_{\Omega}\rho_\varepsilon^a|\nabla u|^2\right)^{1/2}.$$
We will denote by $\tilde H^1(\Omega,\rho_\varepsilon^a(y)\mathrm{d}z)$ the closure of $C^\infty_c(\overline\Omega\setminus\Sigma)$ with respect to the norm $\|\cdot\|_{H^1(\Omega,\rho_\varepsilon^a(y)\mathrm{d}z)}$. In particular, when $a<1$, there is a natural isometry (on balls $B$ centered in a point on $\Sigma$ of any radius)
$$T_\varepsilon^a: \tilde H^{1}(B,\rho_\varepsilon^a(y)\mathrm{d}z)\to\tilde H^{1}(B): u\mapsto v=\sqrt{\rho_\varepsilon^a} u,$$
where  $\tilde H^{1}(B)$ is endowed with the equivalent norm with squared expression
\[
Q_\varepsilon(v)=\int_B|\nabla v|^2 +\left[ \left(\dfrac{\partial_y \rho_\varepsilon^a}{2\rho_\varepsilon^a}\right)^2+\partial_y\left(\dfrac{\partial_y \rho_\varepsilon^a}{2\rho_\varepsilon^a}\right)\right]v^2-\int_{	\partial B}\dfrac{\partial_y \rho_\varepsilon^a}{2\rho_\varepsilon^a}yv^2\;,
\]
(this is detailed in the appendix \ref{app:isometries}). We remark that both in the super singular and super degenerate cases, that is when $a\in(-\infty,-1]\cup[1,+\infty)$ and $\varepsilon=0$, when the weight is taken outside the $A_2$-Muckenhoupt class,  one has 
\begin{equation}\label{tildeH1=H1}
H^{1,a}(\Omega)=\tilde H^{1,a}(\Omega)\,.
\end{equation}
This happens for very opposite reasons: roughly speaking, when $a\leq-1$ then the singularity is so strong to force the function to annihiliate on $\Sigma$ (we will call this case the super singular case). Instead, when $a\geq1$, then the strong degeneracy leaves enough freedom to the function to allow it to be very irregular through $\Sigma$ (we will call this case the super degenerate case). In the latter case, $\Sigma$ has vanishing capacity with respect to the energy $\int |y|^a |\nabla u|^2$.

The Sobolev embedding theorems are stated in details in \cite{SirTerVit1} as inequalities which are uniform in $\varepsilon$. This point is fundamental in order to develop a local regularity theory which is stable with respect to the regularization parameter  $\varepsilon$. Hence, following some results contained in \cite{Haj}, the critical Sobolev exponents do depend on how the weighted measures $\mathrm{d}\mu=\rho_\varepsilon^a(y)\mathrm{d}z$ scale on balls of small radius $r>0$: one can check that there exists $b,d>0$ independent from $\varepsilon\geq0$ (in the locally integrable case $a>-1$) such that for small radii
$$\mu(B_r(z))\geq br^d.$$
So, we can define the effective dimension
$$d=n+1+a^+=n^*(a),$$
and the Sobolev optimal exponent is
$$2^*(a)=\frac{2d}{d-2}=\frac{2(n+1+a^+)}{n+a^+-1}.$$
For details one can refer to Theorems 2.4 and 2.5 in \cite{SirTerVit1}.

In the very same way one can define weighted Sobolev spaces for the class of weights $\omega_\varepsilon^a$; that is, the spaces $H^1(\Omega,\omega_\varepsilon^a(y)\mathrm{d}z)=\tilde H^1(\Omega,\omega_\varepsilon^a(y)\mathrm{d}z)$ (the equality is due to the fact that $\omega_\varepsilon^a$ is always a super degenerate weight as $a<1$) and $H^1_0(\Omega,\omega_\varepsilon^a(y)\mathrm{d}z)$.

In this case one can check that there exist two positive constants $\overline b,\overline d>0$ independent on $\varepsilon\geq0$ such that $\mathrm{d}\overline\mu=\omega_\varepsilon^a(y)\mathrm{d}z$ has the following growth condition on small balls of radius $r>0$
$$\overline\mu(B_r(z))\geq \overline br^{\overline d},$$
and the effective dimension is given by
$\overline d=n+1+2+(-a)^+=n+3+(-a)^+=\overline n^*(a)$. Hence one can state the following
\begin{Theorem}\label{sobemb1}
Let $a\in(-\infty,1)$, $n\geq1$, $\varepsilon\geq0$ and $u\in C^1_c(\Omega)$. Then there exists a constant which does not depend on $\varepsilon\geq0$ such that
\begin{equation*}\label{soboBH}
\left(\int_{\Omega}\omega_\varepsilon^a|u|^{\overline2^*(a)}\right)^{2/\overline2^*(a)}\leq c(\overline d,\overline b,\Omega)\int_{\Omega}\omega_\varepsilon^a|\nabla u|^2,
\end{equation*}
where the optimal embedding exponent is
\begin{equation*}\label{2*a}
\overline2^*(a)=\frac{2\overline d}{\overline d-2}=\frac{2(n+3+(-a)^+)}{n+(-a)^++1}.
\end{equation*}
\end{Theorem}

\subsection{Energy solutions}
Throughout the paper, we are going to consider different elliptic equations depending on different families of weights. Nevertheless, we will deal with right hand sides having forcing terms, terms expressed by the divergence of a given field and drift terms (we will see that any other possible term that will appear can be translated in one of these). In order to give an unified definition of energy solutions to weighted problems, we will consider a generic measurable weight function $w$, and define an energy solution $u$ in $B_1$ to
\begin{equation}\label{divP}
-\mathrm{div}\left(w A\nabla u\right)=wf+\mathrm{div}\left(wF\right)+w\,b\cdot\nabla u\qquad\mathrm{in \ }B_1.
\end{equation}
We say that $u\in H^1(B_1,w\mathrm{d}z)$ is an energy solution to \eqref{divP} if
\begin{equation}\label{variationP}
\int_{B_1}wA(x,y)\nabla u\cdot\nabla\phi=\int_{B_1}wf\phi-\int_{B_1} wF\cdot\nabla\phi+\int_{B_1}w(b\cdot\nabla u)\phi,\qquad\forall\phi\in C^{\infty}_c(B_1)\cap H^1(B_1,w\mathrm{d}z),
\end{equation}
any time the terms in the right hand side give sense to the previous integrals. We remark that we are not assuming local integrability of the weight, and this is the reason why we must consider test functions in the suitable weighted Sobolev space.

Now, we recall the consequent definition of energy solutions in case the weight term is given by $\rho_\varepsilon^a(y)$, with $a\in\R$ and $\varepsilon\geq0$ (the following definition is contained in \cite{SirTerVit1}). Let us consider the following problem
\begin{equation}\label{La}
-\mathrm{div}\left(\rho_\varepsilon^a A\nabla u\right)=\rho_\varepsilon^af+\mathrm{div}\left(\rho_\varepsilon^aF\right)\qquad\mathrm{in \ }B_1.
\end{equation}
We say that $u\in H^{1}(B_1,\rho_\varepsilon^a(y)\mathrm{d}z)$ is an energy solution to \eqref{La} if
\begin{equation}\label{variationLa}
\int_{B_1}\rho_\varepsilon^aA(x,y)\nabla u\cdot\nabla\phi=\int_{B_1}\rho_\varepsilon^af\phi-\int_{B_1} \rho_\varepsilon^aF\cdot\nabla\phi,\qquad\forall\phi\in C^{\infty}_c(B_1)\cap H^{1}(B_1,\rho_\varepsilon^a(y)\mathrm{d}z).
\end{equation}
We remark that the condition in \eqref{variationLa} can be equivalently expressed testing with any $\phi\in C^\infty_c(B_1\setminus\Sigma)$ if $a\in(-\infty,-1]\cup[1,+\infty)$ and $\varepsilon=0$.
In order to give a sense to energy solutions to \eqref{La} we need the following minimal hypothesis on the right hand side.
\begin{Assumption}[H$f\rho_\varepsilon^a$]
Let $a\in(-1,+\infty)$. Then if $n\geq2$ or $n=1$ and $a^+>0$, the forcing term $f$ in \eqref{La} belongs to $L^p(B_1,\rho_\varepsilon^a(y)\mathrm{d}z)$ with $p\geq(2^*(a))'$ the conjugate exponent of $2^*(a)$; that is,
$$(2^*(a))'=\frac{2(n+1+a^+)}{n+a^++3}.$$
If $n=1$ and $a^+=0$ then $f\in L^p(B_1,\rho_\varepsilon^a(y)\mathrm{d}z)$ with $p>1$.

Let $a\in(-\infty,-1]$. Then if $n\geq2$, the condition on the forcing term is $(\rho_\varepsilon^a)^{1/2}f\in L^p(B_1)$ with $p\geq(2^*(a))'=(2^*)'$. If $n=1$, then any $p>1$ is allowed.
\end{Assumption}
\begin{Assumption}[H$F\rho_\varepsilon^a$]
Let $a\in(-1,+\infty)$. The condition on the field $F=(F^1,...,F^{n+1})$ in \eqref{La} is $F\in L^p(B_1,\rho_\varepsilon^a(y)\mathrm{d}z)$ with $p\geq2$. Let $a\in(-\infty,-1]$. Then the condition is $(\rho_\varepsilon^a)^{1/2}F\in L^p(B_1)$ with $p\geq2$.
\end{Assumption}
\subsection{Some preliminary results on the auxiliary equation}

We are concerned with local regularity of energy odd solutions to \eqref{La} with $a\in(-\infty,1)$ and $\varepsilon\geq0$. Our analysis relies in the validity of suitable Liouville type theorems which hold true whenever $a>-1$; that is, when the weight $|y|^a$ is locally integrable. In order to ensure regularity results also in the super singular case $a\leq-1$, we will consider the ratio $w$ between the odd solution $u$ and the function $v^a_\varepsilon$ defined in \eqref{veps}
which is odd and satisfies
\begin{equation}\label{vepsa}
\mathrm{div}\left(\rho_\varepsilon^aA\cdot\nabla v_\varepsilon^a\right)=\mathrm{div}_x\left(\rho_\varepsilon^a\mu\tilde B\cdot \nabla_xv_\varepsilon^a\right)+\mathrm{div}_x(T) \qquad\mathrm{in \ } B_1,
\end{equation}
whenever the right hand side in the equation satisfies suitable integrability assumptions and the matrix $A$ is written as in Notation (HA+). As we have already remarked in the introduction, such a function $v^a_\varepsilon$ plays the role of the characteristic odd solution $y|y|^{-a}$ in presence of a matrix and of regularization.

The following Lemma is a formal computation
\begin{Lemma}\label{rhofrac}
Let $a\in\mathbb{R}$, $\varepsilon>0$ and let $u,v$ be solutions to
$$-\mathrm{div}\left(\rho_\varepsilon^aA\nabla u\right)=\rho_\varepsilon^af,\qquad-\mathrm{div}\left(\rho_\varepsilon^aA\nabla v\right)=\rho_\varepsilon^ag\qquad\mathrm{in \ } B_1,$$
with $v>0$ and $A$ satisfying Assumption (HA). Then the function $w=u/v$ is solution to
\begin{equation*}
-\mathrm{div}\left(\rho_\varepsilon^av^2A\nabla w\right)=\rho_\varepsilon^avf-\rho_\varepsilon^aug\qquad\mathrm{in \ } B_1.
\end{equation*}
\end{Lemma}
\proof
Let recall $\rho=\rho_\varepsilon^a$. Then
\begin{eqnarray*}
-\mathrm{div}\left(\rho v^2A\nabla w\right)&=&-\mathrm{div}\left(\rho v^2A\left(\frac{\nabla u}{v}-\frac{u\nabla v}{v^2}\right)\right)\nonumber\\
&=&-\mathrm{div}\left(\rho vA\nabla u-\rho uA\nabla v\right)\nonumber\\
&=&-v\mathrm{div}\left(\rho A\nabla u\right)-\rho\nabla v\cdot(A\nabla u)+u\mathrm{div}\left(\rho A\nabla v\right)+\rho\nabla u\cdot(A\nabla v)\nonumber\\
&=&-v\mathrm{div}\left(\rho A\nabla u\right)-\rho\nabla v\cdot(A\nabla u)+u\mathrm{div}\left(\rho A\nabla v\right)+\rho\nabla v\cdot(A^T\nabla u)\nonumber\\
&=&\rho vf-\rho ug.
\end{eqnarray*}
\endproof

The new class of weights appearing in the auxiliary equation for the ratio $w=u/v^a_\varepsilon$ is given by $\rho_\varepsilon^a(v_\varepsilon^a)^2$ and it will be equivalent (using \eqref{boundmu}) to
\begin{equation*}
\omega_\varepsilon^a(y)=\rho_\varepsilon^a(y)\pi_\varepsilon^a(y)=\rho_\varepsilon^a(y)(1-a)^2(\chi_\varepsilon^a(y))^2=\rho_\varepsilon^a(y)\left((1-a)\int_0^y\rho_\varepsilon^{-a}(s)\mathrm{d}s\right)^2.
\end{equation*}
We remark that, considering $a\in(-\infty,1)$, such a class of weights is always super degenerate; that is, at $\Sigma$
$$\omega_\varepsilon^a(y)\sim\begin{cases}
y^2 &\mathrm{if \ }\varepsilon>0\\
|y|^{2-a} &\mathrm{if \ }\varepsilon=0,
\end{cases}$$
with $2-a\in(1,+\infty)$.

Formal computations show that the auxiliary equation for $w$ (which corresponds to equation \eqref{BHLaw} in Proposition \ref{prop1} for $\varepsilon=0$ and $A=\mathbb I$) in $B_r$ for any $r<1$ is given by 
\begin{equation}\label{BHepseq}
-\mathrm{div}\left(\rho_\varepsilon^a(v_\varepsilon^a)^2A\nabla w\right)=\rho_\varepsilon^a(v_\varepsilon^a)^2\left(\overline f +Vw-\frac{\overline F\cdot\nabla v_\varepsilon^a}{v_\varepsilon^a}\right)+\mathrm{div}\left(\rho_\varepsilon^a(v_\varepsilon^a)^2\overline F\right)\qquad\mathrm{in \ } B_r,
\end{equation}
with
$$\overline f:=\frac{f}{v_\varepsilon^a},\qquad\overline F:=\frac{F}{v_\varepsilon^a}$$
and
$$V:=\frac{\mathrm{div}_x\left(\mu\tilde B\cdot \nabla_xv_\varepsilon^a\right)}{v_\varepsilon^a}+\frac{\mathrm{div}_x(T)}{\rho_\varepsilon^av_\varepsilon^a}.$$
Actually we can rewrite the $0$-order term, obtaining that the auxiliary equation for $w$ in $B_r$ is given by 
\begin{eqnarray}\label{BHepseq1}
-\mathrm{div}\left(\rho_\varepsilon^a(v_\varepsilon^a)^2A\nabla w\right)&=&\rho_\varepsilon^a(v_\varepsilon^a)^2\left(\overline f -\frac{\overline F\cdot\nabla v_\varepsilon^a}{v_\varepsilon^a}\right)+\mathrm{div}\left(\rho_\varepsilon^a(v_\varepsilon^a)^2\overline F\right)\nonumber\\
&&+\mathrm{div}_x\left(\rho_\varepsilon^a(v_\varepsilon^a)^2b^{\tilde A}w\right)-\rho_\varepsilon^a(v_\varepsilon^a)^2\left(b^{\tilde A}\cdot b^{\mathbb I}w+b^{\tilde A}\cdot\nabla_xw\right)\nonumber\\
&&+\mathrm{div}_x\left(\rho_\varepsilon^a(v_\varepsilon^a)^2\overline Tw\right)-\rho_\varepsilon^a(v_\varepsilon^a)^2\left(\overline T\cdot b^{\mathbb I}w+\overline T\cdot\nabla_xw\right),
\end{eqnarray}
where for a $(n,n)$-dimensional matrix $M$
\begin{equation*}
b^M=M\cdot\frac{\nabla_xv_\varepsilon^a}{v_\varepsilon^a},\qquad\mathrm{and}\qquad \overline T=\frac{T}{\rho_\varepsilon^av_\varepsilon^a}.
\end{equation*}
Thus we can write the equation the following form:
\begin{eqnarray}\label{BHepseq2}
-\mathrm{div}\left(\rho_\varepsilon^a(v_\varepsilon^a)^2A\nabla w\right)&=&\rho_\varepsilon^a(v_\varepsilon^a)^2f+\mathrm{div}\left(\rho_\varepsilon^a(v_\varepsilon^a)^2 F_1\right)\nonumber\\
&&+\mathrm{div}_x\left(\rho_\varepsilon^a(v_\varepsilon^a)^2F_2w\right)+\rho_\varepsilon^a(v_\varepsilon^a)^2Vw+\rho_\varepsilon^a(v_\varepsilon^a)^2b\cdot\nabla_xw.
\end{eqnarray}
We would like to prove that $w$ is an even energy solution to \eqref{BHepseq2} in $B_r$ in the sense that $w\in H^1(B_r,\omega_\varepsilon^a(y)\mathrm{d}z)$ and satisfies
\begin{eqnarray*}
\int_{B_r}\rho_\varepsilon^a(v_\varepsilon^a)^2A\nabla w\cdot\nabla\phi&=&\int_{B_r}\rho_\varepsilon^a(v_\varepsilon^a)^2f\phi-\int_{B_r}\rho_\varepsilon^a(v_\varepsilon^a)^2 F_1\cdot\nabla\phi\nonumber\\
&&-\int_{B_r}\rho_\varepsilon^a(v_\varepsilon^a)^2F_2w\cdot\nabla_x\phi\nonumber\\
&&+\int_{B_r}\rho_\varepsilon^a(v_\varepsilon^a)^2Vw\phi+\int_{B_r}\rho_\varepsilon^a(v_\varepsilon^a)^2(b\cdot\nabla_xw)\phi,
\end{eqnarray*}
for any $\phi\in C^{\infty}_c(B_r\setminus\Sigma)$ (as we have already remarked, super degeneracy allows us to take test functions compactly supported away from $\Sigma$). In order to give a sense to energy solutions to \eqref{BHepseq2} we need the following minimal hypothesis on the right hand side.

\begin{Assumption}[H$f\omega_\varepsilon^a$]
Let $a\in(-\infty,1)$. Then the forcing term $f$ in \eqref{BHepseq2} belongs to $L^p(B_1,\omega_\varepsilon^a(y)\mathrm{d}z)$ with $p\geq(\overline 2^*(a))'$ the conjugate exponent of $\overline 2^*(a)$; that is,
$$(\overline 2^*(a))'=\frac{2(n+3+(-a)^+)}{n+(-a)^++5}.$$
\end{Assumption}
\begin{Assumption}[H$F_1\omega_\varepsilon^a$]
Let $a\in(-\infty,1)$. Then the field term $F_1$ in \eqref{BHepseq2} belongs to $L^p(B_1,\omega_\varepsilon^a(y)\mathrm{d}z)$ with $p\geq2$.
\end{Assumption}
\begin{Assumption}[H$F_2\omega_\varepsilon^a$]
Let $a\in(-\infty,1)$. Then the field term $F_2$ in \eqref{BHepseq2} belongs to $L^p(B_1,\omega_\varepsilon^a(y)\mathrm{d}z)$ with $p\geq \overline d=n+3+(-a)^+$.
\end{Assumption}
\begin{Assumption}[H$V\omega_\varepsilon^a$]
Let $a\in(-\infty,1)$. Then the $0$-order term $V$ in \eqref{BHepseq2} belongs to $L^p(B_1,\omega_\varepsilon^a(y)\mathrm{d}z)$ with
$$p\geq\frac{\overline d}{2}=\frac{n+3+(-a)^+}{2}.$$
\end{Assumption}
\begin{Assumption}[H$b\omega_\varepsilon^a$]
Let $a\in(-\infty,1)$. Then the field $b$ the drift term in \eqref{BHepseq2} belongs to $L^p(B_1,\omega_\varepsilon^a(y)\mathrm{d}z)$ with $p\geq \overline d=n+3+(-a)^+$.
\end{Assumption}

We will need the following important result, which contains also Proposition \ref{prop1} when $\varepsilon=0$ and $A=\mathbb{I}$.

\begin{Proposition}\label{BoundaryHarnackeps}
Let $a\in(-\infty,1)$, $\varepsilon\geq0$ and let $u_\varepsilon\in H^1(B_1,\rho_\varepsilon^a(y)\mathrm{d}z)$ be an odd energy solution to \eqref{La} in $B_1$. Then, fixed $0<r<1$, the function $w_\varepsilon=u_\varepsilon/v_\varepsilon^a$ is an even energy solution in $H^1(B_r,\omega_\varepsilon^a(y)\mathrm{d}z)$ to \eqref{BHepseq1}, provided that the right hand side satisfies the suitable integrability assumptions stated above.
\end{Proposition}
\proof
First, we want to show that $w_\varepsilon\in H^1(B_r,\rho_\varepsilon^a(y)(v_\varepsilon^a)^2(x,y)\mathrm{d}z)$. We remark that since $\frac{1}{C}\leq\mu\leq C$ and since the weight is super degenerate, we have that at $\Sigma$
$$\rho_\varepsilon^a(v_\varepsilon^a)^2\sim\omega_\varepsilon^a\sim\begin{cases}
|y|^{2-a} & \mathrm{if \ }\varepsilon=0\\
|y|^2 & \mathrm{if \ }\varepsilon>0,
\end{cases}$$
with $2-a\in(1,+\infty)$, then the (H=W)-property does not necessarily hold (due to the lack of a Poincar\'e inequality, see \cite{SirTerVit1}). Nevertheless, we can argue as follows: let $\eta\in C^\infty_c(B_1)$ be a radial decreasing cut off function such that $0\leq\eta\leq1$ and $\eta\equiv1$ in $B_r$. Let also for $\delta>0$
$$f_\delta(y)=\begin{cases}
0 & \mathrm{in \ }B_1\cap\{|y|\leq\delta\} \\
\log\frac{y}{\delta} &  \mathrm{in \ }B_1\cap\{\delta\leq|y|\leq\delta e\} \\
1 &  \mathrm{in \ }B_1\cap\{\delta e\leq|y|\}.
\end{cases}$$
Let $\varphi_\delta=\eta f_\delta$, then $|\varphi_\delta|\leq1$ and $|\nabla\varphi_\delta|\leq c/y$ uniformly in $\delta>0$. We remark that one can replace $f_\delta$ with a function with the same properties which is $C^\infty(B_1)$. So,
\begin{equation}\label{L2delta}
\int_{B_1}\rho_\varepsilon^a(v_\varepsilon^a)^2|\varphi_\delta w_\varepsilon|^2\leq\int_{B_1}\rho_\varepsilon^au_\varepsilon^2\leq c.
\end{equation}
Obviously in $B_1\setminus\Sigma$ equation \eqref{BHepseq1} holds.
It is an easy consequence of Lemma \ref{rhofrac}, using that $v_\varepsilon^a$ is an odd energy solution to \eqref{vepsa} in $B_1$ and that $v_\varepsilon^a>0$ in $B_1\setminus\Sigma$. Then, testing the equation with $\varphi_\delta^2w_\varepsilon$, we obtain
\begin{eqnarray}\label{H1delta}
\int_{B_1}\rho_\varepsilon^a(v_\varepsilon^a)^2A\nabla(\varphi_\delta w_\varepsilon)\cdot\nabla(\varphi_\delta w_\varepsilon)&=&\int_{B_1}\rho_\varepsilon^a(v_\varepsilon^a)^2\left(\varphi_\delta^2A\nabla w_\varepsilon\cdot\nabla w_\varepsilon+2\varphi_\delta w_\varepsilon A\nabla w_\varepsilon\cdot\nabla\varphi_\delta+w_\varepsilon^2A\nabla\varphi_\delta\cdot\nabla\varphi_\delta\right)\nonumber\\
&=&\int_{B_1}(RHS)\varphi_\delta^2w_\varepsilon+\int_{B_1}\rho_\varepsilon^a(v_\varepsilon^a)^2w_\varepsilon^2A\nabla\varphi_\delta\cdot\nabla\varphi_\delta\nonumber\\
&\leq&\int_{B_1}(RHS)\varphi_\delta^2w_\varepsilon+c\int_{B_1}\frac{\rho_\varepsilon^a}{y^2}u_\varepsilon^2\leq c,
\end{eqnarray}
and this is true by the weighted Hardy inequality in \eqref{hard}, weighted Sobolev embeddings (Theorem \ref{sobemb1}) and Assumptions $(H f\omega_\varepsilon^a)$, $(H F_1\omega_\varepsilon^a)$, $(H F_2\omega_\varepsilon^a)$, $(H V\omega_\varepsilon^a)$ and $(H b\omega_\varepsilon^a)$ which give a bound on the term with $(RHS)$ of equation \eqref{BHepseq1}. We remark that, fixed $\delta>0$, the boundedness in norm $H^1(B_1,\rho_\varepsilon^a(y)(v_\varepsilon^a)^2(x,y)\mathrm{d}z)$ is enough to ensure that $\varphi_\delta w_\varepsilon$ belongs to the same space. In fact, they have compact support away from $\Sigma$, and hence these norms are equivalent to the usual $H^1$-norm. Since the bounds in \eqref{L2delta} and \eqref{H1delta} are uniform in $\delta>0$, this is enough to have weak convergence for the sequence $\varphi_\delta w_\varepsilon$ in $H^1(B_r,\rho_\varepsilon^a(y)(v_\varepsilon^a)^2(x,y)\mathrm{d}z)$ as $\delta\to0$ and of course the limit is $w_\varepsilon$ (it is almost everywhere pointwise limit).\\\\
We have already remarked that in $B_r\setminus\Sigma$ equation \eqref{BHepseq1} holds. Then, one can conclude since the weighted Sobolev space $H^1(B_r,\rho_\varepsilon^a(y)(v_\varepsilon^a)^2(x,y)\mathrm{d}z)$ is super degenerate, and consequently test functions can be taken in $C^\infty_c(B_r\setminus\Sigma)$.
\endproof

\section{Liouville type theorems}\label{sec:liouville}
In this section we present two important results which will be our main tool in order to prove regularity local estimates which are uniform with respect to $\varepsilon\geq0$.
\begin{Theorem}\label{Liouvilleodd}
Let $a\in(-1,1)$, $\varepsilon\geq 0$ and $w$ be a solution to
\begin{equation*}
\begin{cases}
-\mathrm{div}(\rho_\varepsilon^a(y)\nabla w)=0 & \mathrm{in \ }\mathbb{R}^{n+1}_+\\
w(x,0)=0,
\end{cases}
\end{equation*}
and let us suppose that for some $\gamma\in[0,1-a)$, $C>0$ it holds
\begin{equation}\label{eq:wronggrowth0}
|w(z)|\leq C (1 + |z|^{\gamma})
\end{equation}
for every $z$. Then $w$ is identically zero.
\end{Theorem}
\proof
It is enough to prove the result only for $\varepsilon\in\{0,1\}$. In fact for any other value of $\varepsilon>0$ we can normalize the problem falling in the case $\varepsilon=1$.\\\\
$\bf{Case \ 1:}$ \ $\varepsilon=0$.\\
Let us consider $w\in H^{1,a}_{\mathrm{loc}}(\mathbb{R}^{n+1}_+)$ satisfying the conditions of the statement, that is, solution in the following sense
$$\int_{\mathbb{R}^{n+1}_+}y^a\nabla w\cdot\nabla\phi=0\qquad\forall\phi\in C^\infty_c(\mathbb{R}^{n+1}_+).$$
Let us define
$$E(r)=\frac{1}{r^{n+a-1}}\int_{B_r^+}y^{a}|\nabla w|^2,\qquad H(r)=\frac{1}{r^{n+a}}\int_{\partial^+ B_r^+}y^aw^2.$$
Note that, as the weight $y^a$ is locally integrable, \eqref{eq:wronggrowth0} implies
\begin{equation}\label{eq:wgh0}
H(r)\leq C(1+r^{2\gamma})\;, \forall r>0\;.
\end{equation}
Now, defining $w^r(x)=w(rx)$,  we have
$$E(r)=\int_{B_1^+}y^a|\nabla w^r|^2\qquad\mathrm{and}\qquad H(r)=\int_{S^n_+}y^a(w^r)^2,$$
and hence 
$$H'(r)=\frac{2}{r}E(r).$$
We are looking for the best constant in the following trace Poincar\'e inequality
\begin{equation}\label{tpti}
\int_{B_1^+}y^a|\nabla u|^2\geq\lambda(a)\int_{S^n_+}y^a u^2.
\end{equation}
Actually we are able to provide the best constant $\lambda(a)$ in \eqref{tpti}, since $u(x,y)=y^{1-a}$ is the unique function in $\tilde H^{1,a}(B_1^+)$ which solves
\begin{equation*}
\begin{cases}
-L_au=0 &\mathrm{in} \ B_1^+\\
u>0 &\mathrm{in} \ B_1^+\\
u(x,0)=0\\
\nabla u\cdot\nu=\lambda(a)u &\mathrm{in} \ S_+^{n},
\end{cases}
\end{equation*}
with $\lambda(a)=1-a$. However $\lambda(a)$ is the same of \eqref{lama}. Hence $H'(r)\geq \frac{2\lambda(a)}{r}H(r)$, and integrating, there we infer
\begin{equation*}
\frac{H(r)}{r^{2(1-a)}}\geq H(1),
\end{equation*}
we obtain that if $w$ is not trivial, its growth at infinity is at least $r^{1-a}$, in contradiction with \eqref{eq:wgh0} taking $r$ large.\\
$\bf{Case \ 2:}$ \ $\varepsilon=1$.\\
Let us define
$$E(r)=\frac{1}{r^{n+a-1}}\int_{B_r^+}(1+y^2)^{a/2}|\nabla w|^2,\qquad H(r)=\frac{1}{r^{n+a}}\int_{\partial^+ B_r^+}(1+y^2)^{a/2}w^2.$$
Note that, as $a>-1$, the $\rho_\varepsilon^{a}$'s are uniformly locally integrable and thus \ref{eq:wronggrowth0} implies again
\begin{equation}\label{eq:wgh05}
H(r)\leq C(1+r^{2\gamma})\;, \forall r>0\;, \text{with $\gamma<1-a$.}
\end{equation}
Hence,
\begin{equation}\label{H11}
H'(r)=\frac{2}{r}E(r)-\frac{a}{r^{n+a+1}}\int_{\partial^+ B_r^+}(1+y^2)^{a/2-1}w^2.
\end{equation}
Moreover, defining $w^r(x)=w(rx)$ one has
$$E(r)=\int_{B_1^+}\left(\frac{1}{r^2}+y^2\right)^{a/2}|\nabla w^r|^2\qquad\mathrm{and}\qquad H(r)=\int_{S^n_+}\left(\frac{1}{r^2}+y^2\right)^{a/2}(w^r)^2.$$
By Lemma \ref{A1} and Remark \ref{A2}, one can find for any radius $r>0$ the best constant $\lambda_r(a)$ such that
\begin{equation}\label{tptir}
\int_{B_1^+}\left(\frac{1}{r^2}+y^2\right)^{a/2}|\nabla u|^2\geq\lambda_r(a)\int_{S^n_+}\left(\frac{1}{r^2}+y^2\right)^{a/2} u^2.
\end{equation}
Defining $\rho_k(y)=\left(\frac{1}{r_k^2}+y^2\right)^{a/2}$ with $r_k\to+\infty$ as $k\to+\infty$, one can see
$$\lambda(a)=\min_{v\in\tilde H^1(B_1^+)\setminus\{0\}}\frac{Q_a(v)}{\int_{S^n_+}v^2}\qquad\mathrm{and}\qquad\lambda_k(a)=\min_{v\in\tilde H^1(B_1^+)\setminus\{0\}}\frac{Q_{\rho_k}(v)}{\int_{S^n_+}v^2}.$$
By Lemma \ref{A1}, $\lambda_k(a)\to\lambda(a)=1-a$ as $k\to+\infty$.\\\\
Now we want to prove that the correction term in \eqref{H11} is of lower order as $r\to+\infty$. By \eqref{tracehard}, we have that in $\tilde C^\infty_c(B_1^+)$
\begin{equation*}
\int_{B_1^+}\rho_r|\nabla u|^2\geq c_0\int_{\partial B_1^+}\frac{\rho_r}{y}u^2.
\end{equation*}
Hence
\begin{eqnarray*}
\left|\frac{a}{r^{n+a+1}}\int_{\partial^+ B_r^+}\left(1+y^2\right)^{a/2-1}w^2\right|&\leq&\frac{|a|}{r^{n+a+1}}\int_{\partial^+ B_r^+}\left(1+y^2\right)^{a/2-1/2}w^2\\
&=&\frac{|a|}{r^{2}}\int_{S^n_+}\left(\frac{1}{r^2}+y^2\right)^{a/2-1/2}(w^r)^2\\
&\leq&\frac{|a|}{r^{2}}\int_{S^n_+}\left(\frac{1}{r^2}+y^2\right)^{a/2}y^{-1}(w^r)^2\\
&\leq&\frac{|a|}{c_0r^{2}}\int_{B_1^+}\left(\frac{1}{r^2}+y^2\right)^{a/2}|\nabla w^r|^2\\
&=&\frac{|a|}{c_0r^{2}}E(r).
\end{eqnarray*}
Hence for $r$ large enough
\begin{equation*}
H'(r)\geq \frac{2\lambda_r(a)}{r}H(r),
\end{equation*}
and since $\lambda_r(a)\to\lambda(a)=1-a$, by integrating the above expression we deduce that, for all small $\delta$, there exists $r_0>0$ such that, for every $r>r_0$
\begin{equation*}
\frac{H(r)}{r^{2(1-a-\delta)}}\geq H(r_0),
\end{equation*}
which says that if $w$ is not trivial, its growth at infinity is at least $r^{1-a-\delta}$. Taking $\delta>0$ so small that $1-a-\delta>\gamma$ we find a contradiction with \eqref{eq:wgh05}.
\endproof

\begin{Theorem}\label{LiouvilleBH}
Let $a\in(-\infty,1)$, $\varepsilon\geq0$ and $w$ be a solution to
\begin{equation*}
\begin{cases}
-\mathrm{div}((\omega_\varepsilon^a(y))^{-1}\nabla w)=0 & \mathrm{in \ }\mathbb{R}^{n+1}_+\\
w=0 &\mathrm{in \ }\mathbb{R}^{n}\times\{0\},
\end{cases}
\end{equation*}
and let us suppose that for some $\gamma\in[0,1)$, $C>0$ it holds
\begin{equation}\label{eq:wronggrowthBH}
|w(z)|\leq C\omega_\varepsilon^a(y)(1 + |z|^\gamma)
\end{equation}
for every $z=(x,y)$. Then $w$ is identically zero.
\end{Theorem}
\proof
By a simple normalization argument, it is enough to prove the result only for $\varepsilon\in\{0,1\}$. We start with\\\\
$\bf{Case \ 1:}$ \ $\varepsilon=0$.\\
The case falls into the proof of $\bf{Case \ 1}$ in \cite[Theorem 3.4]{SirTerVit1} replacing $a\in(-\infty,1)$ with $(a-2)\in(-\infty,-1)$.\\\\
$\bf{Case \ 2:}$ \ $\varepsilon=1$.\\
Let us now define
$$E(r)=\frac{1}{r^{n+(a-2)-1}}\int_{B_r^+}(\omega_1^a(y))^{-1}|\nabla w|^2,\qquad\textrm{and}\qquad H(r)=\frac{1}{r^{n+(a-2)}}\int_{\partial^+ B_r^+}(\omega_1^a(y))^{-1}w^2.$$
Note that, defining $w^r(x)=w(rx)$ one has
$$E(r)=\int_{B_1^+}(\omega_{1/r}^a(y))^{-1}|\nabla w^r|^2,\qquad\textrm{and}\qquad H(r)=\int_{S^n_+}(\omega_{1/r}^a(y))^{-1}(w^r)^2.$$
First we remark that the growth condition \eqref{eq:wronggrowthBH} implies the following upper bound
\begin{equation}\label{eq:wgh}
H(r)\leq Cr^{-2(a-2)}(1+r^{2\gamma})\;, \forall r>0\;,
\end{equation}
(due to the local integrability of  $y^{2-a}$) and heence 
$$\int_{S^n_+}\omega_{1/r}^a(y)\leq c$$
uniformly in $r>0$. Therefore,
\begin{equation}\label{H111}
H'(r)=\frac{2}{r}E(r)+\int_{S^n_+}\frac{d}{dr}[(\omega_{1/r}^a(y))^{-1}](w^r)^2.
\end{equation}
By Lemma \ref{B1} and Lemma \ref{A2BH}, one can find for any radius $r>0$ the best constant $\mu_r(a)$ such that
\begin{equation}\label{tptir2}
\int_{B_1^+}(\omega_{1/r}^a(y))^{-1}|\nabla u|^2\geq\mu_r(a)\int_{S^n_+}(\omega_{1/r}^a(y))^{-1}|u|^2.
\end{equation}
Defining $(\omega^a_k(y))^{-1}=(\omega_{1/r_k}^a(y))^{-1}$ and $\mu_k=\mu_{r_k}$ with $r_k\to+\infty$ as $k\to+\infty$, by Lemma \ref{A2BH}, $\mu_k(a)\to\mu(a)=1-(a-2)=3-a$ as $k\to+\infty$.\\\\
Now we want to prove that the correction term in \eqref{H111} is of lower order as $r\to+\infty$. By \eqref{tracehardBH}, we have that\begin{equation*}
\int_{B_1^+}(\omega_{1/r}^a(y))^{-1}|\nabla u|^2\geq c_0\int_{S^n_+}\frac{(\omega_{1/r}^a(y))^{-1}}{y}u^2.
\end{equation*}
Using
\begin{eqnarray*}
\int_{S^n_+}\frac{d}{dr}[(\omega_{1/r}^a(y))^{-1}](w^r)^2&=&\frac{a}{r^3}\int_{S^n_+}\frac{1}{\left(\frac{1}{r^2}+y^2\right)}(\omega_{1/r}^a(y))^{-1}(w^r)^2\\
&&-\frac{2a}{r^3}\int_{S^n_+}\frac{\int_0^y\left(\frac{1}{r^2}+s^2\right)^{-a/2-1}}{\int_0^y\left(\frac{1}{r^2}+s^2\right)^{-a/2}}(\omega_{1/r}^a(y))^{-1}(w^r)^2,
\end{eqnarray*}
we can estimate the first term of the rest as follows
\begin{eqnarray*}
\left|\frac{a}{r^3}\int_{S^n_+}\frac{1}{\left(\frac{1}{r^2}+y^2\right)}(\omega_{1/r}^a(y))^{-1}(w^r)^2\right|&\leq&\frac{|a|}{r^2}\int_{S^n_+}\frac{1}{\left(\frac{1}{r^2}+y^2\right)^{1/2}}(\omega_{1/r}^a(y))^{-1}(w^r)^2\\
&\leq&\frac{|a|}{r^2}\int_{S^n_+}\frac{(\omega_{1/r}^a(y))^{-1}}{y}(w^r)^2\leq\frac{c}{r^{2}}E(r).
\end{eqnarray*}
Moreover, when $a\leq0$ the second term of the rest can be estimated as
\begin{eqnarray*}
\left|\frac{2a}{r^3}\int_{S^n_+}\frac{\int_0^y\left(\frac{1}{r^2}+s^2\right)^{-a/2-1}}{\int_0^y\left(\frac{1}{r^2}+s^2\right)^{-a/2}}(\omega_{1/r}^a(y))^{-1}(w^r)^2\right|&\leq&\frac{2|a|}{r^2}\int_{S^n_+}\frac{ry\int_0^{ry}\left(1+s^2\right)^{-a/2-1}}{\int_0^{ry}\left(1+s^2\right)^{-a/2}}\frac{(\omega_{1/r}^a(y))^{-1}}{y}(w^r)^2\\
&\leq&\frac{2|a|}{r^2}\int_{S^n_+}\frac{(\omega_{1/r}^a(y))^{-1}}{y}(w^r)^2\leq\frac{c}{r^{2}}E(r),
\end{eqnarray*}
and this is due to the fact that, calling $z=ry\in[0,+\infty)$, by the fact that
$$f(z)=\frac{z\int_0^{z}\left(1+s^2\right)^{-a/2-1}}{\int_0^{z}\left(1+s^2\right)^{-a/2}}$$
is continuous and such that $f(0)=0$ and
$$f(z)\sim_{z\to+\infty}\begin{cases}cz^a &\mathrm{if \ }a\in(-1,0]\\ \frac{\log z}{z} &\mathrm{if \ }a=-1\\ \frac{1}{z} &\mathrm{if \ }a<-1\end{cases}$$
and hence $f(z)\leq c$ in $[0,+\infty)$. Instead, when $a\in(0,1)$ the second term of the rest can be estimated as
\begin{eqnarray*}
\left|\frac{2a}{r^3}\int_{S^n_+}\frac{\int_0^y\left(\frac{1}{r^2}+s^2\right)^{-a/2-1}}{\int_0^y\left(\frac{1}{r^2}+s^2\right)^{-a/2}}(\omega_{1/r}^a(y))^{-1}(w^r)^2\right|&\leq&\frac{2|a|}{r^{2-a}}\int_{S^n_+}\frac{(ry)^{1-a}\int_0^{ry}\left(1+s^2\right)^{-a/2-1}}{\int_0^{ry}\left(1+s^2\right)^{-a/2}}\frac{(\omega_{1/r}^a(y))^{-1}}{y^{1-a}}(w^r)^2\\
&\leq&\frac{2|a|}{r^{2-a}}\int_{S^n_+}\frac{(\omega_{1/r}^a(y))^{-1}}{y^{1-a}}(w^r)^2\leq\frac{2|a|}{r^{2-a}}\int_{S^n_+}\frac{(\omega_{1/r}^a(y))^{-1}}{y}(w^r)^2\\
&\leq&\frac{c}{r^{2-a}}E(r),
\end{eqnarray*}
using the fact that $0\leq y\leq1$ and by the fact that
$$f(z)=\frac{z^{1-a}\int_0^{z}\left(1+s^2\right)^{-a/2-1}}{\int_0^{z}\left(1+s^2\right)^{-a/2}}$$
is continuous and such that $f(0)=0$ and
$$f(z)\sim_{z\to+\infty}c$$
and hence $f(z)\leq c$ in $[0,+\infty)$.
Hence for $r$ large enough
\begin{equation*}
H'(r)\geq \frac{2\mu_r(a)}{r}H(r),
\end{equation*}
and since $\mu_r(a)\to\mu(a)=1-(a-2)$, we can choose a small $\delta>0$ such that $1-(a-2)-\delta>\gamma-(a-2)$. Hence, by integrating the above expression we deduce that there exists $r_0>0$ such that, for every $r>r_0$, we have $\mu_r(a)>1-(a-2)-\delta>\gamma-(a-2)$ and
\begin{equation*}
\frac{H(r)}{r^{2(1-(a-2)-\delta)}}\geq H(r_0),
\end{equation*}
which is in contradiction with \eqref{eq:wgh} for $r$ large if $w$ is not trivial.
\endproof

\begin{Corollary}\label{LiouvilleBHbis}
Let $a\in(-\infty,1)$, $\varepsilon\geq0$ and $w$ be a solution to
\begin{equation*}
\begin{cases}
-\mathrm{div}(\omega_\varepsilon^a(y)\nabla w)=0 & \mathrm{in \ }\mathbb{R}^{n+1}_+\\
\omega_\varepsilon^a\partial_yw=0 &\mathrm{in \ }\mathbb{R}^{n}\times\{0\},
\end{cases}
\end{equation*}
and let us suppose that for some $\gamma\in[0,1)$, $C>0$ it holds
\begin{equation}\label{gro0BH}
|w(z)|\leq C(1 + |z|^\gamma)
\end{equation}
for every $z$. Then $w$ is constant.
\end{Corollary}
\proof
Again, it is enough to treat the cases $\varepsilon\in\{0,1\}$. Let us assume $\varepsilon=1$, the case $\varepsilon=0$ coincides with the case $\varepsilon=0$ in \cite[Corollary 3.5]{SirTerVit1}, by replacing in the proof $a\in(-1,+\infty)$ by $(2-a)\in(1,+\infty)$. Then we have (by an even reflection across $\Sigma$) an even solution $w$  to
\begin{equation*}
-\mathrm{div}\left(\omega_1^a(y)\nabla w\right)=0 \qquad\mathrm{in} \ \mathbb{R}^{n+1}.
\end{equation*}
Such a solution is $w\in H^{1,2}_{\mathrm{loc}}(\mathbb{R}^{n+1})=H^{1}_{\mathrm{loc}}(\mathbb{R}^{n+1},|y|^2\mathrm{d}z)$, with the growth condition \eqref{gro0BH}. Now we observe that, as $w$ is not constant with a sublinear growth at infinity,  $v=\omega_1^a(y)\partial_y w$ can not be trivial, otherwise
$w$ would be globally harmonic and sublinear, in contradiction with the Liouville theorem in \cite{NorTavTerVer}. Hence, if $w$ is not constant, $v$  must be  an odd and nontrivial solution to
\begin{equation*}
\begin{cases}
-\mathrm{div}\left((\omega_1^a(y))^{-1}\nabla v\right)=0 &\mathrm{in \ }\mathbb{R}^{n+1}_+\\
v=0 &\mathrm{in} \ \{y=0\}.
\end{cases}
\end{equation*}
By the arguments in the proof of  Theorem \ref{LiouvilleBH}, we know that the weighted average of $v^2$ must satisfy a minimal growth rate as
\begin{equation*}\label{gro1}
 H(r)=\frac{1}{r^{n+(a-2)}}\int_{\partial^+ B_r^+}(\omega_1^a(y))^{-1}v^2\geq cr^{2(1-(a-2)-\delta)}, \qquad 1-\delta>\gamma\;,
\end{equation*}
for $r\geq r_0$ depending on $\delta>0$ chosen. Therefore, by integrating, we obtain
\[
\int_{B_r^+}\omega_1^a(y)(\partial_y w)^2=\int_0^r\mathrm{d}t\int_{\partial^+ B_t^+}(\omega_1^a(y))^{-1}v^2\geq c r^{n-(a-2)+2-2\delta}\;.
\]
On the other hand, we have, by \eqref{gro0BH}
\begin{equation*}
\begin{split}
\int_{B_r^+}\omega_1^a(y) (\partial_y w)^2\leq \int_{B_r^+}\omega_1^a(y) |\nabla w|^2\\ \leq c  \int_{B_{2r}^+}\omega_1^a(y) |w|^2\leq c(1+r^{n-(a-2)+2\gamma})
\end{split}
\end{equation*}
in contradiction with the previous inequality when $r$ is large, since $1-\delta>\gamma$.
\endproof

\section{Local uniform bounds in H\"older spaces for the auxiliary problem}
As a first step in our regularity theory for odd solutions, we state some results on local uniform estimates for solutions to \eqref{simpleBH}; that is,
\begin{equation*}
-\mathrm{div}\left(\rho_\varepsilon^a(v_\varepsilon^a)^2A\nabla u_\varepsilon\right)=\rho_\varepsilon^a(v_\varepsilon^a)^2f_\varepsilon+\mathrm{div}\left(\rho_\varepsilon^a(v_\varepsilon^a)^2F_\varepsilon\right)+\rho_\varepsilon^a(v_\varepsilon^a)^2b_\varepsilon\cdot\nabla u_\varepsilon\qquad\mathrm{in \ }B_1.
\end{equation*}
Using a Moser iteration argument (see also \cite[Section 8.4]{GilTru}), one can prove the following nowadays standard result.
\begin{Proposition}\label{Moser}
Let $a\in(-\infty,1)$ and $\varepsilon\geq0$. Let $u\in H^{1}(B_1,\omega_\varepsilon^a(y)\mathrm{d}z)$ be an energy solution to \eqref{simpleBH}. Let $\beta>1$,
$$p_1>\frac{\overline d}{2}=\frac{n+3+(-a)^+}{2}, \qquad p_2,p_3>\overline d.$$
Let moreover $\|b\|_{L^{p_3}(B_1,\omega_\varepsilon^a(y)\mathrm{d}z)}\leq b$.
Then, for any $0<r<1$ there exists a positive constant independent of $\varepsilon$ (depending on $n$, $a$, $p_1$, $p_2$, $p_3$, $\beta$, $b$, $r$ and $\alpha$) such that
$$\|u\|_{L^\infty(B_{r})}\leq c\left(\|u\|_{L^{\beta}(B_1,\omega_\varepsilon^a(y)\mathrm{d}z)}+\|f\|_{L^{p_1}(B_1,\omega_\varepsilon^a(y)\mathrm{d}z)}+\|F\|_{L^{p_2}(B_1,\omega_\varepsilon^a(y)\mathrm{d}z)}\right).$$
\end{Proposition}
\proof
The proof follows the same steps as in \cite[Proposition 2.17]{SirTerVit1}, but iterating the Sobolev embedding in Theorem \ref{sobemb1}.
\endproof

\subsection{Local uniform bounds in $C^{0,\alpha}$ spaces}
\begin{Theorem}\label{C0alphaBH}
Let $a\in(-\infty,1)$ and as $\varepsilon\to0$ let $\{u_\varepsilon\}$ be a family of solutions in $B_1^+$ of \eqref{simpleBH} satisfying boundary conditions (evenness)
$$
\rho_\varepsilon^a(v_\varepsilon^a)^2\partial_yu_\varepsilon=0 \qquad \mathrm{on \ }\partial^0 B_1^+.
$$
Let $r\in(0,1)$, $\beta>1$, $p_1>\frac{\overline d}{2}=\frac{n+3+(-a)^+}{2}$, $p_2,p_3>\overline d$, and $\alpha\in(0,1)\cap(0,2-\frac{n+3+(-a)^+}{p_1}]\cap(0,1-\frac{n+3+(-a)^+}{p_2}]\cap(0,1-\frac{n+3+(-a)^+}{p_3}]$. Let $\|b\|_{L^{p_3}(B_1,\omega_\varepsilon^a(y)\mathrm{d}z)}\leq b$. Let moreover $A$ satisfy assumption (HA) with continuous coefficients. There is a positive constant depending on $a$, $n$, $\beta$, $p_1$, $p_2$, $p_3$, $b$, $\alpha$ and $r$ only such that
$$\|u_\varepsilon\|_{C^{0,\alpha}(B_r^+)}\leq c\left(\|u_\varepsilon\|_{L^\beta(B_1^+,\omega_\varepsilon^a(y)\mathrm{d}z)}+ \|f_\varepsilon\|_{L^{p_1}(B_1^+,\omega_\varepsilon^a(y)\mathrm{d}z)}+\|F_\varepsilon\|_{L^{p_2}(B_1,\omega_\varepsilon^a(y)\mathrm{d}z)}\right).$$
\end{Theorem}
\proof
The proof follows the very same steps as in the proof of \cite[Theorem 4.1]{SirTerVit1}.
First, one has to remark that the suitable H\"older continuity for $\varepsilon\geq0$ fixed is given by the theory developed for even solutions to degenerate problems in \cite{SirTerVit1}. Then, one can argue by contradiction with the usual blow up argument considering two blow up sequences
\begin{equation*}
v_k(z)=\frac{(\eta u_k)(z_k+r_kz)-(\eta u_k)(z_k)}{L_kr_k^\alpha},\qquad w_k(z)=\frac{\eta(z_k)(u_k(z_k+r_kz)- u_k(z_k))}{L_kr_k^\alpha},
\end{equation*}
(with the same asymptotic behaviour on compact sets) defined in the rescaled domains $B(k)=\frac{B-z_k}{r_k}$ (where $B=B_{\frac{1+r}{2}}$ and $\{z_k\}$ is one of the two sequences of points where H\"older seminorms blow up), the first possessing some uniform H\"older continuity, and the second one satisfying suitable problems on rescaled domains which blow up. In order to complete the proof we prove some steps.\\\\
{\bf Step 1: blow-ups.} The first thing to do is to characterize the possible asymptotic behaviours of the weights $\rho_\varepsilon^a(v_\varepsilon^a)^2$ in the rescaled points: that is,
\begin{eqnarray*}
p_k(z)&:=&\rho_\varepsilon^a(y_k+r_ky)(v_\varepsilon^a(z_k+r_kz))^2\nonumber\\
&=&\left(\varepsilon_k^2+(y_k+r_ky)^2\right)^{a/2}\left(\int_0^{y_k+r_ky}\left(\varepsilon_k^2+s^2\right)^{-a/2}\mu(x_k+r_kx,s)^{-1}\mathrm{d}s\right)^2.
\end{eqnarray*}
To this end, let us define by $\Gamma_k=(\varepsilon_k,y_k,r_k)$ and denote  $\nu_k=|\Gamma_k|$. The latter is a bounded sequence and, up to subsequences, has finite limit $\nu=|(0,y_\infty,0)|\geq0$ (where we have assumed $z_k\to z_\infty=(x_\infty,y_\infty)$). Taking possibly another subsequence, we may assume that  the normalized sequence
$$\tilde\Gamma_k=\frac{\Gamma_k}{\nu_k}=(\tilde\varepsilon_k,\tilde y_k,\tilde r_k)=\left(\frac{\varepsilon_k}{\nu_k},\frac{y_k}{\nu_k},\frac{r_k}{\nu_k}\right)$$
has a limit
$$\tilde\Gamma_k\to\tilde\Gamma=(\tilde\varepsilon,\tilde y,\tilde r)\in S^2\;,$$
and moreover that
$$ \lim_{k\to+\infty}\frac{\tilde y_k}{\tilde r_k}=\tilde l\in[0,+\infty].$$
Thus we can consider $\tilde\Sigma=\lim\Sigma_k$; that is,
\begin{equation*}
\tilde\Sigma=\begin{cases}
\{(x,y)\in\mathbb{R}^{n+1} \ : \ y=-\tilde l \} &\mathrm{if}\; \tilde l <+\infty,\\
\emptyset &\mathrm{if}\; \tilde l =+\infty.
\end{cases}
\end{equation*}
After rescaling the independent variables, we find new weights having the form:
\begin{equation*}
p_k(z)=\nu_k^{2-a}\left(\tilde\varepsilon_k^2+\left(\tilde y_k+\tilde r_ky\right)^2\right)^{a/2}\left(\int_0^{\tilde y_k+\tilde r_ky}\left(\tilde\varepsilon_k^2+t^2\right)^{-a/2}\mu(x_k+r_kx,\nu_kt)^{-1}\mathrm{d}t\right)^2,
\end{equation*}
and, in order to study their asymptotics, we have to distinguish between different cases:\\\\
{\bf Case 1.} $\nu>0$. Then, $\tilde r=\tilde\varepsilon=0$ and $\tilde y=1$. Moreover, it is easy to see, using that $1/\mu$ is continuous, that $p_k(z)= c+o(1)$.\\\\
{\bf Case 2.} $\nu=0$ and $\tilde\varepsilon=0$ ($\tilde y\neq0\lor\tilde r\neq0$). Using the continuity of $1/\mu$, up to a vertical translation of $-\tilde l$, one obtains
$$p_k(z)=\nu_k^{2-a}\tilde p(y)(1+o(1))$$
where
$$\tilde p(y)=\begin{cases}
c &\mathrm{if \ }\tilde r=0,\\
c|y|^{2-a} &\mathrm{if \ }\tilde r\neq 0.
\end{cases}$$
{\bf Case 3.} $\nu=0$ and $\tilde\varepsilon\in(0,1)$. Using again the continuity of $1/\mu$, , up to a vertical translation of $-\tilde l$, we obtain
$$p_k(z)=\nu_k^{2-a}\tilde p(y)(1+o(1))$$
where
$$\tilde p(y)=\begin{cases}
c &\mathrm{if \ }\tilde r=0,\\
c \ \omega_1^a(y) &\mathrm{if \ }\tilde r\neq 0.
\end{cases}$$
in the second case up to a dilation of $\frac{\tilde\varepsilon}{\tilde r}$.\\\\
{\bf Case 4.} $\nu=0$ and $\tilde\varepsilon=1$ ($\tilde y=0\land\tilde r=0$). As usual, by the continuity of $1/\mu$, , up to a vertical translation of $-\tilde l$, one obtains, if $\tilde r_k=o(\tilde y_k)$
$$p_k(z)=\nu_k^{2-a}\tilde y_k^2c(1+o(1)),$$
and otherwise
$$p_k(z)=\nu_k^{2-a}\tilde r_k^2c|y|^2(1+o(1)).$$
\\\\
Let us define
$$h_k=\begin{cases}
\nu_k^{2-a} &\mathrm{in \ {\bf Cases \ 1,2,3}},\\
\nu_k^{2-a}\tilde y_k^2 &\mathrm{in \ {\bf Case \ 4}, \ and \ } \tilde r_k=o(\tilde y_k)\\
\nu_k^{2-a}\tilde r_k^2 &\mathrm{in \ {\bf Case \ 4}, \ otherwise,}
\end{cases}$$
and $\tilde p_k=\frac{p_k}{h_k}$. We have shown that, up to the suitable normalization, the rescaled weights $\tilde p_k$ do converge uniformly  to $\tilde p$ on compact sets of $\R^{n+1}\setminus\tilde\Sigma$ (or the whole $\R^{n+1}$ whenever $\tilde\Sigma=\emptyset$). Note that this latter case is equivalent to the limiting weight $\tilde p$ be constant. \\\\
{\bf Step 2: the limiting equation and uniform-in-$k$ energy estimates.}
The equation for the rescaled variable $w_k$ becomes:
\begin{eqnarray}\label{eq:rescaled}
-\mathrm{div}(\tilde p_kA(z_k+r_k \cdot)\nabla w_k)(z)&=& \frac{\eta(z_k)}{L_k}r_k^{2-\alpha}\tilde p_k(z)f_k(z_k+r_kz)\nonumber\\
&&+\mathrm{div}\left(\frac{\eta(z_k)}{L_k}r_k^{1-\alpha} \tilde p_k(\cdot)F_k(z_k+r_k\cdot)\right)(z)\nonumber\\
&&+r_k\tilde p_k(z)b_k(z_k+r_kz)\cdot\nabla w_k(z).
\end{eqnarray}
By a Caccioppoli type inequality, easily obtained by multiplying \eqref{eq:rescaled} by $\overline\eta^2 w_k$,  being $\overline\eta$ a cut-off function, taking into account that the $w_k$ are uniformly bounded and that 
\begin{itemize}
\item the first term in the right hand side of \eqref{eq:rescaled} is bounded in $L^1_{loc}$;
\item the  field $\frac{\eta(z_k)}{L_k}r_k^{1-\alpha} F_k(z_k+r_k\cdot)$ in the second term in the right hand side of \eqref{eq:rescaled} is bounded in $L^2_{loc}(\tilde p(z)\mathrm{d}z)$;
\item $r_kb_k(z_k+r_k\cdot)$ is bounded in $L^2_{loc}(\tilde p_k(z)\mathrm{d}z)$;
\end{itemize}
then we obtain uniform-in-$k$ energy bounds holding on compact subsets of $\R^{n+1}$:
\[
\forall R>0,\;\exists c>0,\;\forall k, \qquad
\int_{B_R}\tilde p_k A(z_k+r_k z) \nabla w_k\cdot \nabla w_k\leq c\;.
\]
The computations are very similar to the ones done in the following step.\\\\
{\bf Step 3: the right hand side vanishes as $k\to+\infty$.}
Next we wixh to check that the right hand sides in the rescaled equations vanish in $L^1_{\mathrm{loc}}(\R^{n+1}\setminus\tilde\Sigma)$ (or $L^1_{\mathrm{loc}}(\R^{n+1})$ whenever $\tilde\Sigma=\emptyset$), and that consequently the limit $w$ is an energy solution of
\begin{equation}\label{limittilde}
-\mathrm{div}\left(\tilde pA(z_\infty)\nabla w\right)=0\qquad \mathrm{in} \ \mathbb{R}^{n+1}\setminus\tilde\Sigma,
\end{equation}
even with respect to $\tilde\Sigma$ (when not empty). We use the continuity of the matrix $A$ in order to obtain a constant coefficients matrix in the limit equation \eqref{limittilde} together with the fact that $\tilde\Sigma$ is invariant with respect to the limit matrix to have evenness across the characteristic hyperplane. \\\\
Let us show that the right hand sides vanish in $L^1_{\mathrm{loc}}$ at least for one of the cases (the other cases are very similar), for instance when $h_k=\nu_k^{2-a}\tilde y_k^2$; that is, {\bf Case 4}, when $\tilde r_k=o(\tilde y_k)$. Indeed, let $\phi\in C^\infty_c(\mathbb{R}^{n+1})$: using the fact that for $k$ large enough $\mathrm{supp}(\phi)\subset B_R\subset B(k)$, using H\"older inequality, we have
\begin{eqnarray}
&&\left|\int_{B_R}\tilde p_k(z) f_{\varepsilon_k}(z_k+r_kz)\phi(z)\mathrm{d}z\right|\nonumber\\
&\leq&\|\phi\|_{L^\infty(B_R)}\left(\frac{1}{r_k^{n+1}h_k}\int_{B_{r_kR}(z_k)}\rho_{\varepsilon_k}^a(\zeta_{n+1})(v_{\varepsilon_k}^a(\zeta))^2|f_{\varepsilon_k}(\zeta)|^{p_1}\mathrm{d}\zeta\right)^{1/p_1}\left(\int_{B_{R}}\tilde p_k(z)\mathrm{d}z\right)^{1/p'_1}\nonumber\\
&\leq&cr_k^{-\frac{n+1}{p_1}}\nu_k^{-\frac{2-a}{p_1}}\tilde y_k^{-\frac{2}{p_1}},\nonumber
\end{eqnarray}
and hence the first term in the right hand side converges to zero since $\alpha\leq2-\frac{n+3+(-a)^+}{p_1}$, $\tilde r_k=\frac{r_k}{\nu_k}$, the fact that $0\leq r_k\leq \nu_k$ and having
$$\frac{\eta(z_k)}{L_k}r_k^{2-\alpha-\frac{n+3+(-a)^+}{p_1}}\left(\frac{r_k^{(-a)^+}}{\nu_k^{-a}}\right)^{1/p_1}\left(\frac{\tilde r_k}{\tilde y_k}\right)^{2/p_1}\to0.$$
With analogous computations one can check that also the second term in the right hand side vanishes.\\\\
Concerning the third term, one can estimate the integral as follows
\begin{eqnarray}\label{3p}
&&\\
&&r_k\left|\int_{B_R}\tilde p_k(z)b_k(z_k+r_kz)\cdot\nabla w_k(z)\phi(z)\right|\nonumber\\
&\leq& r_k\|\phi\|_{L^\infty(B_R)}\left(\int_{B_R}\tilde p_k|\nabla w_k|^2\right)^{\frac{1}{2}}\left(\int_{B_R}\tilde p_k\right)^{\frac{p_3-2}{2p_3}}\left(\frac{1}{h_kr_k^{n+1}}\int_{B_{r_kR}(z_k)} \rho_{\varepsilon_k}^a(\zeta_{n+1})(v_{\varepsilon_k}^a(\zeta))^2|b_{\varepsilon_k}(\zeta)|^{p_3}\right)^{\frac{1}{p_3}}\nonumber\\
&\leq&cb^{\frac{1}{p_3}}r_k^{1-\frac{n+3+(-a)^+}{p_3}}\left(\frac{r_k^{(-a)^+}}{\nu_k^{-a}}\right)^{1/p_3}\left(\frac{\tilde r_k}{\tilde y_k}\right)^{2/p_3}\left(\int_{B_R}\tilde p_k|\nabla w_k|^2\right)^{1/2}\nonumber\\
&=&t_k\left(\int_{B_R}\tilde p_k|\nabla w_k|^2\right)^{1/2}\nonumber.
\end{eqnarray}
The sequence $t_k$ converges to zero, having $p_3>n+3+(-a)^+$. Moreover, the full term vanishes using the uniform energy bound obtained in the previous step.\\\\
{\bf Step 4: the limit belongs to $H^1_{loc}(\R^{n+1},\tilde p\mathrm{d}z)$.}  At this point, always up to subsequences, we know that the (pointwise) convergence to $w$ holds also in the weak $H^1_{loc}(\R^{n+1}\setminus\tilde\Sigma)$ topology. Now we wish to  infer that the limit $w$ belongs to the space $H^1_{loc}(\R^{n+1},\tilde p \mathrm{d}z)$ as the closure of $C^\infty$ with respect to the weighted norm (as defined in \S\ref{sect:sobolev}).  Let us start with the easiest case when  $\tilde\Sigma=\emptyset$ and the limiting weight $\tilde p$ is constant. Moreover,  we know tha $\tilde p_k$ converge uniformly to $\tilde p$ on compact sets. Thus the sequence $w_k$ converges weakly $H^1$ to $w$ on each compact subset. The convergence to The case when $\tilde \Sigma\neq \emptyset$ requires a more thorough analysis.
%
%
In order to ensure that $w\in H^1_{\mathrm{loc}}(\R^{n+1},\tilde p(y)\mathrm{d}z)$ also when $\tilde\Sigma\neq \emptyset$, one can argue as follows: using the fact that $\mu$ is continuous with $\frac{1}{C}<\mu<C$, then fixed a compact set of $\R^{n+1}$, we can find positive constants $c_k,C_k$ (which are uniformly bounded from above and below by two constants respectively $0<c_1<c_2<+\infty$) such that 
$$c_k\tilde\omega_{\varepsilon_k}(y_k+r_ky)\leq \tilde p_k(x,y)\leq C_k\tilde \omega_{\varepsilon_k}(y_k+r_ky)\qquad\mathrm{and}\qquad\frac{C_k}{c_k}\to1.$$
where we have denoted 
\begin{equation}\label{eq:tildeomega}
\tilde\omega_k=\frac{\omega_{\varepsilon_k}}{h_k}.
\end{equation} 
Now, reabsorbing the weights as in \eqref{Qomega1+} and using the family of isometries given by
$$\overline T_k(w_k)=\left(\tilde\omega_{k}(y_k+r_ky)\right)^{1/2}w_k=W_k,$$
one obtains uniform boundedness of the $W_k$'s  in $ H^1_{\mathrm{loc}}(\R^{n+1})$, and hence they weakly converge in the same space to $W$. Coming back with the inverse isometry associated with the limit weight
$$\overline T(w)=(\tilde p(y))^{1/2}w=W,$$
we obtain $w\in H^1_{\mathrm{loc}}(\R^{n+1},\tilde p(y)\mathrm{d}z)$.\\\\
{\bf{Step 5: end of the proof.}} Next we wish to show that $w$ solves the equation in \eqref{limittilde} also across the limiting characteristic hyperplane $\tilde \Sigma$.  Indeed, using $w\in H^1_{\mathrm{loc}}(\R^{n+1},\tilde p(y)\mathrm{d}z)$ jointly with equation \eqref{limittilde} holding  in $\R^{n+1}\setminus\Sigma$, using that  $C^\infty_c(\overline B_R\setminus \Sigma)$ is actually dense in $H^1(B_R,\tilde p(y)\mathrm{d}z)$ (all the weights here, including the limit one, are super degenerate), as we have already remarked in \eqref{tildeH1=H1}.
 Eventually, one can reach a contradiction by applying the suitable Liouville theorems in \cite{NorTavTerVer, SirTerVit1} and Corollary \ref{LiouvilleBHbis} for the case $\tilde p(y)=\omega_1^a(y)$.
\endproof

\subsection{Local uniform bounds in $C^{1,\alpha}$ spaces}
\begin{Theorem}\label{C1alphaBH}
Let $a\in(-\infty,1)$ and as $\varepsilon\to0$ let $\{u_\varepsilon\}$ be a family of solutions in $B_1^+$ of \eqref{simpleBH} satisfying boundary conditions (evenness)
$$
\rho_\varepsilon^a(v_\varepsilon^a)^2\partial_yu_\varepsilon=0 \qquad \mathrm{on \ }\partial^0 B_1^+.
$$
Let $r\in(0,1)$, $\beta>1$, $p_1,p_2>\overline d=n+3+(-a)^+$. Let $\|b\|_{L^{2p_2}(B_1,\omega_\varepsilon^a(y)\mathrm{d}z)}\leq b$. Let $F_\varepsilon=(F^1_\varepsilon,...,F^{n+1}_\varepsilon)$ with the $y$-component vanishing on $\Sigma$: $F^{n+1}_\varepsilon(x,0)= F^y_\varepsilon(x,0)=0$ in $\partial^0B_1^+$. Let moreover $A$ satisfy assumption (HA) with $\alpha$-H\"older continuous coefficients and $\alpha\in(0,1-\frac{n+3+(-a)^+}{p_1}]\cap(0,1-\frac{n+3+(-a)^+}{p_2}]$. There is a positive constant depending on $a$, $n$, $\beta$, $p_1$, $p_2$, $b$, $\alpha$ and $r$ only such that
$$\|u_\varepsilon\|_{C^{1,\alpha}(B_r^+)}\leq c\left(\|u_\varepsilon\|_{L^\beta(B_1^+,\omega_\varepsilon^a(y)\mathrm{d}z)}+ \|f_\varepsilon\|_{L^{p_1}(B_1^+,\omega_\varepsilon^a(y)\mathrm{d}z)}+\|F_\varepsilon\|_{C^{0,\alpha}(B_1^+)}\right).$$
\end{Theorem}
\proof
We wish to follow the same steps of proof of \cite[Theorems 5.1 and 5.2]{SirTerVit1}. Among others, we have to deal with an additional difficulty; that is, our weights here do depend on the full variable $z=(x,y)$ and not on $y$ only. For our purposes, we can take advantage of the fact that the ratio
\begin{equation}
\gamma_\varepsilon^a(x,y)=:\frac{v_\varepsilon^a(x,y)}{\chi_\varepsilon^a(y)}=\frac{\int_0^y\rho_\varepsilon^{-a}(s)\mu^{-1}(x,s)\mathrm{d}s}{\int_0^y\rho_\varepsilon^{-a}(s)\mathrm{d}s}
\end{equation}
is uniformly bounded in $C^{0,\alpha}$ with respect to $\varepsilon$ (just apply Lemma \ref{c0alphaG}, using the fact that $\mu^{-1}\in C^{0,\alpha}$ since the matrix $A$ possesses $\alpha$-H\"older continuous coefficients). Hence, one can rewrite our operator as
\begin{equation}
\textrm{div}(\rho_\varepsilon^a(y)(v_\varepsilon^a(x,y))^2A(x,y)\nabla u_\varepsilon)=\textrm{div}(\omega_\varepsilon^a(y)A_\varepsilon(x,y)\nabla u_\varepsilon),
\end{equation}
where, up to constants, the new family of matrices is defined as
\begin{equation}
A_\varepsilon(x,y)=(\gamma_\varepsilon^a(x,y))^2A(x,y),
\end{equation}
with coefficients which are uniformly bounded in $C^{0,\alpha}$ with respect to $\varepsilon$.\\\\
With these premises, we are now able to follow the construction made in \cite[Theorems 5.1 and 5.2]{SirTerVit1}. Just to give the idea, the contradiction argument uses two blow-up sequences
\begin{equation*}
v_k(z)=\frac{\eta(\hat{z}_k+r_kz)}{L_kr_k^{1+\alpha}}\left(u_k(\hat{z}_k+r_kz)-u_k(\hat{z}_k)\right),\qquad w_k(z)=\frac{\eta(\hat{z}_k)}{L_kr_k^{1+\alpha}}\left(u_k(\hat{z}_k+r_kz)-u_k(\hat{z}_k)\right),
\end{equation*}
for $z\in B(k):=\frac{B-\hat{z}_k}{r_k}$. Hence, one has to work with
$$\overline v_k(z)=v_k(z)-\nabla v_k(0)\cdot z,\qquad\overline w_k(z)=w_k(z)-\nabla w_k(0)\cdot z,$$
or
$$\overline v_k(z)=v_k(z)-\nabla_x v_k(0)\cdot x,\qquad\overline w_k(z)=w_k(z)-\nabla_x w_k(0)\cdot x,$$
respectively when $\frac{d(z_k,\Sigma)}{r_k}\to+\infty$ (in this case we choose $\hat{z}_k=z_k$), or $\frac{d(z_k,\Sigma)}{r_k}\leq c$ uniformly in $k$ (in this case we choose $\hat{z}_k=(x_k,0)$ to be the projection on $\Sigma$ of $z_k$, where $z_k=(x_k,y_k)$).\\\\
Hence, reasoning as in the previous Theorem \ref{C0alphaBH}, one can characterize all possible rescalings of the weights (in facts the possible scalings of weights $p_k$ and $\omega_k$ are the same), and prove that the limit $w$ is an energy entire solution to the suitable limiting problem.\\\\
We remark that in order to show that the limit equation has a constant coefficient matrix one has to reason as in \cite[Remark 5.3]{SirTerVit1}, using the $\alpha$-H\"older continuity of coefficients of the matrix (in this case we will invoke the uniform bounds with respect to $\varepsilon$ in $C^{0,\alpha}$ for the coefficients of $A_{\varepsilon_k}$).\\\\
Nevertheless, we need also to deal with drift terms in the rescaled equations, and we wish to show that they vanish once testing with the suitable test function $\phi$ supported in $B_R$. We assume here that $h_k=\nu_k^{2-a}$ (one of the possible cases). Hence, also in this case we use the fact that we know a priori that the sequence $\{u_k\}$ is uniformly locally bounded in $C^{0,\beta}$ spaces, for any choice of $\beta\in(0,1)$ (follows from Theorem \ref{C0alphaBH}). Reasoning as in \cite[Remark 5.3]{SirTerVit1}, this gives the following energy estimate
\begin{equation}
\int_{B_r} p_k(z)|\nabla u_k|^2\leq\frac{c}{r^{2(1-\beta)}}\int_{B_{2r}}p_k(z).
\end{equation}
Hence, we can estimate
\begin{eqnarray*}
&&r_k\left\vert \int_{B_R}\tilde p_k(z)b_k(\hat z_k+r_kz)\cdot\nabla w_k(z)\phi(z)\right\vert\\
&\leq&\frac{r_k^{1-\alpha}\eta(\hat z_k)}{L_k}\|\phi\|_{L^\infty(B_R)}\int_{B_R}\tilde p_k(z)|b_k(\hat z_k+r_kz)|\cdot|\nabla u_k(\hat z_k+r_kz)|\\
&\leq&\frac{cr_k^{1-\alpha}}{L_k}\left(\int_{B_R}\tilde p_k(z)|b_k(\hat z_k+r_kz)|^2\right)^{1/2}\left(\int_{B_R}\tilde p_k(z)|\nabla u_k(\hat z_k+r_kz)|^2\right)^{1/2}\\
&\leq&\frac{cr_k^{1-\alpha}}{L_k}\left(\int_{B_R}\tilde p_k(z)\right)^{1/2p_2'}\left(\frac{1}{r_k^{n+1}h_k}\int_{B_{r_kR}(\hat z_k)} \rho_{\varepsilon_k}^a(\zeta_{n+1})(v_{\varepsilon_k}^a(\zeta))^2|b_k(\zeta)|^{2p_2}\right)^{1/2p_2}\\
&&\cdot\frac{c}{r_k^{1-\beta}}\left(\int_{2B_R}\tilde p_k(z)\right)^{1/2}\\
&\leq&\frac{c}{L_k}\left(\frac{r_k}{\nu_k}\right)^{1/p_2}\left(\frac{r_k^{(-a)^+}}{\nu_k^{-a}}\right)^{1/2p_2}r_k^{\frac{1}{2}\left(1-\alpha-\frac{n+3+(-a)^+}{p_2}\right)}r_k^{\beta-\frac{1+\alpha}{2}}\to0
\end{eqnarray*}
since $L_k\to+\infty$, $\nu_k\leq r_k$, $p_2>n+3+(-a)^+$, $\alpha\leq1-\frac{n+3+(-a)^+}{p_2}$ and choosing $\beta>\frac{1+\alpha}{2}$.
\endproof

\subsection{Local regularity for the auxiliary equation}
The following is the main result of the paper (we have already stated it in a simplified version in Theorem \ref{holderBHsimple} in the introduction). Let $a\in(-\infty,1)$, the matrix $A$ written as in Notation (HA+) and let $u_\varepsilon$ be an odd energy solution to \eqref{LarhoA} in $B_1^+$; that is,
\begin{equation}\label{1oddBH}
\begin{cases}
-\mathrm{div}\left(\rho_\varepsilon^aA\nabla u_\varepsilon\right)=\rho_\varepsilon^af_\varepsilon+\mathrm{div}\left(\rho_\varepsilon^aF_\varepsilon\right) & \mathrm{in \ } B_1^+\\
u_\varepsilon=0 & \mathrm{on \ }\partial^0B_1^+.
\end{cases}
\end{equation}
Let also $v_\varepsilon^a$ be defined as in \eqref{veps} in $B_1^+$. Then, we have already showed (in Proposition \ref{BoundaryHarnackeps}) that under suitable integrability assumptions on the terms in the right hand side and on coefficients of the matrix $A$, then functions
$$w_\varepsilon=\frac{u_\varepsilon}{v_\varepsilon^a}$$
are even energy solutions (for any $R<1$) to equations \eqref{BHepseq1}; that is,
\begin{eqnarray*}
-\mathrm{div}\left(\rho_\varepsilon^a(v_\varepsilon^a)^2A\nabla w_\varepsilon\right)&=&\rho_\varepsilon^a(v_\varepsilon^a)^2\left(\overline f_\varepsilon -\frac{\overline F_\varepsilon\cdot\nabla v_\varepsilon^a}{v_\varepsilon^a}\right)+\mathrm{div}\left(\rho_\varepsilon^a(v_\varepsilon^a)^2\overline F_\varepsilon\right)\nonumber\\
&&+\mathrm{div}_x\left(\rho_\varepsilon^a(v_\varepsilon^a)^2(b^{\tilde A}_\varepsilon+\overline T_\varepsilon) w_\varepsilon\right)-\rho_\varepsilon^a(v_\varepsilon^a)^2\left((b^{\tilde A}_\varepsilon+\overline T_\varepsilon)\cdot b^{\mathbb I}_\varepsilon w_\varepsilon+(b^{\tilde A}_\varepsilon+\overline T_\varepsilon)\cdot\nabla_xw_\varepsilon\right),
\end{eqnarray*}
with boundary condition
$$\rho_\varepsilon^a(v_\varepsilon^a)^2\partial_yw_\varepsilon=0 \qquad \mathrm{on \ }\partial^0B_R^+,$$
and where we denote by
\begin{equation*}
\overline f_\varepsilon=\frac{f_\varepsilon}{v_\varepsilon^a},\qquad \overline F_\varepsilon=\frac{F_\varepsilon}{v_\varepsilon^a},\qquad b^M_\varepsilon=M\cdot\frac{\nabla_xv_\varepsilon^a}{v_\varepsilon^a},\qquad\mathrm{and}\qquad \overline T_\varepsilon=\frac{T}{\rho_\varepsilon^av_\varepsilon^a},
\end{equation*}
($M$ is a general $(n,n)$-dimensional matrix).

\begin{Theorem}\label{holderBH}
Let $a\in(-\infty,1)$, the matrix $A$ written as in Notation (HA+) and as $\varepsilon\to0$ let $\{u_\varepsilon\}$ be a family of solutions in $B_1^+$ of \eqref{1oddBH}.\\\\
$1)$ Let $r\in(0,1)$, $\beta>1$, $p_1,p_2>\frac{n+3+(-a)^+}{2}$, $p_3,p_4>n+3+(-a)^+$, and $\alpha\in(0,2-\frac{n+3+(-a)^+}{p_1}]\cap(0,2-\frac{n+3+(-a)^+}{p_2}]\cap(0,1-\frac{n+3+(-a)^+}{p_3}]\cap(0,1-\frac{n+3+(-a)^+}{p_4}]$. Let also 
$$\| (b^{\tilde A}_\varepsilon+\overline T_\varepsilon)\cdot b^{\mathbb I}_\varepsilon\|_{L^{p_2}(B_1^+,\omega_\varepsilon^a(y)\mathrm{d}z)}\leq b_1,\qquad \| b^{\tilde A}_\varepsilon+\overline T_\varepsilon\|_{L^{p_4}(B_1^+,\omega_\varepsilon^a(y)\mathrm{d}z)}\leq b_2.$$
Let us moreover take $A$ with continuous coefficients. There is a positive constant depending on $a$, $n$, $\beta$, $p_1$, $p_2$, $p_3$, $p_4$, $b_1$, $b_2$, $\alpha$ and $r$ only such that functions
$$w_\varepsilon=\frac{u_\varepsilon}{v_\varepsilon^a}$$
satisfy
\begin{eqnarray*}
\|w_\varepsilon\|_{C^{0,\alpha}(B_r^+)}&\leq& c\left(\|w_\varepsilon\|_{L^\beta(B_1^+,\omega_\varepsilon^a(y)\mathrm{d}z)}+ \left\|\overline f_\varepsilon-\frac{\overline F_\varepsilon\cdot \nabla v_\varepsilon^a}{v_\varepsilon^a}\right\|_{L^{p_1}(B_1^+,\omega_\varepsilon^a(y)\mathrm{d}z)} +\|\overline F_\varepsilon\|_{L^{p_3}(B_1^+,\omega_\varepsilon^a(y)\mathrm{d}z)}\right).
\end{eqnarray*}
$2)$ Let $r\in(0,1)$, $\beta>1$, $p_1,p_2>n+3+(-a)^+$, and $\alpha\in(0,1-\frac{n+3+(-a)^+}{p_1}]\cap(0,1-\frac{n+3+(-a)^+}{p_2}]$. Let $\overline F_\varepsilon=(\overline F^1_\varepsilon,...,\overline F^{n+1}_\varepsilon)$ with the $y$-component vanishing on $\Sigma$: $\overline F^{n+1}_\varepsilon(x,0)=\overline F^y_\varepsilon(x,0)=0$ in $\partial^0B_1^+$. Let also 
$$\| (b^{\tilde A}_\varepsilon+\overline T_\varepsilon)\cdot b^{\mathbb I}_\varepsilon\|_{L^{p_2}(B_1^+,\omega_\varepsilon^a(y)\mathrm{d}z)}\leq b_1,\qquad \| b^{\tilde A}_\varepsilon+\overline T_\varepsilon\|_{C^{0,\alpha}(B_1^+)}\leq b_2.$$
Let's moreover take $A$ with $\alpha$-H\"older continuous coefficients. There is a positive constant depending on $a$, $n$, $\beta$, $p_1$, $p_2$, $b_1$, $b_2$, $\alpha$ and $r$ only such that functions
$$w_\varepsilon=\frac{u_\varepsilon}{v_\varepsilon^a}$$
satisfy
\begin{eqnarray*}
\|w_\varepsilon\|_{C^{1,\alpha}(B_r^+)}&\leq& c\left(\|w_\varepsilon\|_{L^\beta(B_1^+,\omega_\varepsilon^a(y)\mathrm{d}z)}+\left\|\overline f_\varepsilon-\frac{\overline F_\varepsilon\cdot \nabla v_\varepsilon^a}{v_\varepsilon^a}\right\|_{L^{p_1}(B_1^+,\omega_\varepsilon^a(y)\mathrm{d}z)}+\|\overline F_\varepsilon\|_{C^{0,\alpha}(B_1^+)}\right).
\end{eqnarray*}
\end{Theorem}
We remark that uniform estimates with respect to the regularization are optimal in $C^{1,\alpha}$-spaces (in \cite[Remark 5.4]{SirTerVit1} we provided a counterexample which show that $C^{2,\alpha}$ estimates could not be uniform up to $\Sigma$ as $\varepsilon\to0$).\\\\

In order to prove our main result we have the following useful preliminary result on equations of the form
\begin{equation}\label{BHepseq3}
-\mathrm{div}\left(\rho_\varepsilon^a(v_\varepsilon^a)^2A\nabla u_\varepsilon\right)=\rho_\varepsilon^a(v_\varepsilon^a)^2Vu_\varepsilon+\mathrm{div}\left(\rho_\varepsilon^a(v_\varepsilon^a)^2F u_\varepsilon\right)\qquad\mathrm{in \ } B_1.
\end{equation}
\begin{Lemma}
Let $a\in(-\infty,1)$ and $\varepsilon\geq0$. Let $u\in H^{1}(B_1,\omega_\varepsilon^a(y)\mathrm{d}z)$ be an energy solution to \eqref{BHepseq3}, where $V\in L^{p_1}(B_1,\omega_\varepsilon^a(y)\mathrm{d}z)$ with $p_1>\frac{\overline d}{2}=\frac{n+3+(-a)^+}{2}$ and $F\in L^{p_2}(B_1,\omega_\varepsilon^a(y)\mathrm{d}z)$ with $p_2>\overline d=n+3+(-a)^+$. Let
$$\|V\|_{L^{p_1}(B_1,\omega_\varepsilon^a(y)\mathrm{d}z)}\leq b_1,\qquad \|F\|_{L^{p_2}(B_1,\omega_\varepsilon^a(y)\mathrm{d}z)}\leq b_2.$$
$1)$ Then, for any $0<r<1$ and $\beta>1$ there exists a positive constant independent of $\varepsilon$ (depending on $n$, $a$, $r$, $\beta$, $p_1$, $p_2$, $b_1$, $b_2$), $m_1>\frac{\overline d}{2}$ and $m_2>\overline d$ such that
$$\|Vu\|_{L^{m_1}(B_{r},\omega_\varepsilon^a(y)\mathrm{d}z)}+\|Fu\|_{L^{m_2}(B_{r},\omega_\varepsilon^a(y)\mathrm{d}z)}\leq c\|u\|_{L^{\beta}(B_1,\omega_\varepsilon^a(y)\mathrm{d}z)}.$$
$2)$ If moreover $p_1>\overline d=n+3+(-a)^+$ and $F\in C^{0,\alpha}(B_1)$ for some $\alpha\in(0,1)$, 
$$\|V\|_{L^{p_1}(B_1,\omega_\varepsilon^a(y)\mathrm{d}z)}\leq b_1,\qquad \|F\|_{C^{0,\alpha}(B_1)}\leq b_2,$$
then for any $0<r<1$ and $\beta>1$ there exists a positive constant independent of $\varepsilon$ (depending on $n$, $a$, $r$, $\beta$, $p_1$, $\alpha$, $b_1$, $b_2$), and $m_1>\overline d$ such that
$$\|Vu\|_{L^{m_1}(B_{r},\omega_\varepsilon^a(y)\mathrm{d}z)}+\|Fu\|_{C^{0,\alpha}(B_{r})}\leq c\|u\|_{L^{\beta}(B_1,\omega_\varepsilon^a(y)\mathrm{d}z)}.$$
\end{Lemma}
\proof
The proof is done applying Moser iterations on a finite number of small enough balls which cover $B_r$. The radius of such balls is chosen in order to ensure coercivity of the quadratic forms. Hence, using the fact that the weighted integrability of $V$ and $F$ is suitably large, by a finite number of Moser iterations one can promote the integrability of $u$ itself, up to guarantee that the products $Vu$ and $Fu$ have the desired integrability (this type of argument is classic, see for instance \cite[Section 8.4]{GilTru}). Since the number of iterations is finite, one can control uniformly the constants in the iterative process, proving point 1). At to point 2), thanks to point 1) we can  apply Theorem \ref{C0alphaBH} in order to obtain that the solution is $C^{0,\alpha}$ with a bound which is independent from $\varepsilon$. Hence, we obtain the second inequality taking into account the H\"older continuity of $F$.
\endproof
A relevant consequence of this result is that, under suitable conditions on the $0$-order terms and divergence terms with the solution itself inside (the conditions stated in Theorem \ref{holderBH}), we can treat $Vw$ and $\mathrm{div}(\rho v^2Fw)$ respectively as a fixed forcing term and a divergence term with a given field. As a consequence, we obtain uniform local regularity estimates in Theorem \ref{holderBH} for solutions $w_\varepsilon$ to \eqref{BHepseq} by simply applying Theorems \ref{C0alphaBH} and \ref{C1alphaBH}.

\subsubsection{A criterion for local $C^{1,\alpha}$ estimates}
We would like to show an example of a set of hypothesis for which part $2)$ of our main Theorem \ref{holderBH} holds true; that is, local uniform $C^{1,\alpha}$ estimates for the ratio of odd solutions and the fundamental ones.\\\\
We remark that, as $a\to-\infty$, the decay of the data on $\Sigma$ becomes stronger and stronger.
\begin{Assumption}[$C^{1,\alpha}$]
Let $f_\varepsilon:=y^{\max\{1,1-a\}}g_\varepsilon$ with $g_\varepsilon$ uniformly bounded in $L^{p}(B_1^+,\omega_\varepsilon^a(y)\mathrm{d}z)$ as $\varepsilon\to0$ and
\begin{equation*}
p>n+3+(-a)^+.
\end{equation*}
Let $F_\varepsilon:=y^{\max\{1,1-a\}}G_\varepsilon$ with $G_\varepsilon$ uniformly bounded in $C^{0,\alpha}(B_1^+)$ as $\varepsilon\to0$ and
\begin{equation*}
\alpha>1-\frac{1+\min\{2,2-a\}}{n+3+(-a)^+}.
\end{equation*}
Nevertheless, the matrix $A$, which satisfies Assumption (HA+), must also satisfy some regularity assumptions: $A\in C^{0,\alpha}(B_1^+)$ with $\nabla_x\mu\in C^{0,\alpha}(B_1^+)$ and $T=y\tilde T$ with $\tilde T\in C^{0,\alpha}(B_1^+)$. 
\end{Assumption}

\subsection{Local uniform bounds in H\"older spaces for odd solutions in the $A_2$ case}

Moreover, when the weight is locally integrable; that is, $a\in(-1,1)$, we obtain local estimates for odd solutions working directly on the equation.
\begin{Theorem}\label{holderodd}
Let $a\in(-1,1)$ and as $\varepsilon\to0$ let $\{u_\varepsilon\}$ be a family of solutions in $B_1^+$ of either
\begin{equation}\label{1odd}
-\mathrm{div}\left(\rho_\varepsilon^aA\nabla u_\varepsilon\right)=\rho_\varepsilon^af_\varepsilon+\mathrm{div}\left(\rho_\varepsilon^aF_\varepsilon\right)
\end{equation}
satisfying the Dirichlet boundary condition
\[
u_\varepsilon=0\quad \mathrm{on \ }\partial^0B_1^+.
\]
Let $r\in(0,1)$, $\beta>1$, $p_1>\frac{n+1+a^+}{2}$, $p_2>n+1+a^+$ and $\alpha\in(0,1)\cap(0,1-a)\cap(0,2-\frac{n+1+a^+}{p_1}]\cap(0,1-\frac{n+1+a^+}{p_2}]$. Let moreover $A$ satisfy assumption (HA) with continuous coefficients. There are constants depending on $a$, $n$, $\beta$, $p_1$, $p_2$, $\alpha$ and $r$ only such that
$$\|u_\varepsilon\|_{C^{0,\alpha}(B_r^+)}\leq c\left(\|u_\varepsilon\|_{L^\beta(B_1^+,\rho_\varepsilon^a(y)\mathrm{d}z)}+ \|f_\varepsilon\|_{L^{p_1}(B_1^+,\rho_\varepsilon^a(y)\mathrm{d}z)} + \|F_\varepsilon\|_{L^{p_2}(B_1^+,\rho_\varepsilon^a(y)\mathrm{d}z)}\right).$$
\end{Theorem}
\proof
The proof is obtained by contradiction following the very same passages of \cite[Theorem 4.1]{SirTerVit1}, observing that in presence of the zero Dirichlet boundary condition at $\Sigma$ we obtain a contradiction by applying the Liouville Theorem \ref{Liouvilleodd}. The blow-up sequences invoked are centered in points $\hat z_k\in B^+=B_{\frac{1+r}{2}}\cap\{y\geq0\}$; that is,
\begin{equation*}
v_k(z)=\frac{(\eta u_k)(\hat z_k+r_kz)-(\eta u_k)(\hat z_k)}{L_kr_k^\alpha},\qquad w_k(z)=\frac{\eta(\hat z_k)(u_k(\hat z_k+r_kz)- u_k(\hat z_k))}{L_kr_k^\alpha},
\end{equation*}
with
$$z\in B(k):=\frac{B-\hat z_k}{r_k}.$$
Moreover, if $y_k/r_k\to+\infty$ (where $z_k=(x_k,y_k)$), then we choose $\hat z_k=z_k$, while if $y_k/r_k\leq c$ uniformly, then we choose $\hat z_k=(x_k,0)$. In this second case, we remark that $v_k$ and $w_k$ are antisymmetric with respect to $\{y=0\}$ so that the limit $w$ will be odd in $y$.
\endproof


\appendix

\section{Some special functions}
In this appendix we are going to state and prove some technical results which will allow us to compare, from the regularity point of view, the variable coefficient case with the constant one.
\begin{remark}\label{psia}
Let $a\in(-\infty,1)$, $\varepsilon\geq0$. Then the family of functions
\begin{equation}\label{psi}
\psi_\varepsilon^a(y):=\frac{y\rho_\varepsilon^{-a}(y)}{\int_0^y\rho_\varepsilon^{-a}(s)\mathrm{d}s}
\end{equation}
are monotone in $y$ and uniformly bounded in $L^{\infty}(B_1^+)$ by a constant which does not depend on $\varepsilon$. In fact, denoting $t=y/\varepsilon$, we have

\begin{equation*}
\psi_\varepsilon^a(y)=\psi_1^a\left(\frac{y}{\varepsilon}\right)=\psi_1^a\left(t\right)=\frac{t(1+t^2)^{-a/2}}{\int_0^t(1+s^2)^{-a/2}\mathrm{d}s}.
\end{equation*}
The latter function is continuous and monotone nondecreasing if $a<0$ and nonincreasing if $a\in(0,1)$. Since $\psi_1^a$ has limit $1$ as $t\to0$ and limit $1-a$ as $t\to+\infty$, then
$$\sup_{t>0}\psi_1^a(t)=\max\{1,1-a\}\qquad\mathrm{and}\qquad \inf_{t>0}\psi_1^a(t)=\min\{1,1-a\}.$$

Finally, note that  the family $\psi_\varepsilon^a$ can not be  equicontinuous, nor uniformly bounded in $C^{0,\alpha}(B_1^+)$, while it enjoys the following property:
\begin{equation}\label{eq:property}
\exists c>0\;:\;
\forall \varepsilon\in[0,\varepsilon_0)\;,  \Vert \psi_\varepsilon^a\Vert_{\textrm{Lip}(B_1\cap\{y>\sqrt{\varepsilon}\})}<c\;,
\end{equation}
due to the fact that $\psi_1^a$ is bounded,  has a finite limit as $t\to +\infty$ and  its derivative vanishes as $1/t^2$.

\end{remark}

\begin{Lemma}\label{a2}
Let $a\in(-\infty,1)$, $\varepsilon\geq0$, $\alpha\in(0,1)$ and let $g(x,y,s)\in C^{0,\alpha}_{x,y}(B_1^+)$ uniformly in $s\in[0,y]$, such that $|g(x,y,s)|\leq c|y|^\alpha$ for $(x,y)\in B_1^+$ uniformly in $s\in[0,y]$. Then the family of functions
\begin{equation*}
\mathcal G_\varepsilon(x,y)=\frac{\int_0^y\rho_\varepsilon^{-a}(s)g(x,y,s)\mathrm{d}s}{\int_0^y\rho_\varepsilon^{-a}(s)\mathrm{d}s}
\end{equation*}
is uniformly bounded in $C^{0,\alpha}(B_1^+)$ by a constant which does not depend on $\varepsilon$.
\end{Lemma}
\proof
We remark that the proof follows some ideas of the proof in \cite[Lemma 7.5]{SirTerVit1}, where the case $\varepsilon=0$ is done. The uniform H\"older continuity in the $x$-variable is trivial. Hence, fixed $0<\delta<1$, let us consider the following two sets
$$I_1=\{(y_1,y_2)\ : \ 0\leq y_1\leq y_2<1, \ y_2-y_1\geq\delta y_2\}$$
and
$$I_2=\{(y_1,y_2)\ : \ 0\leq y_1\leq y_2<1, \ y_2-y_1<\delta y_2\}.$$
If we consider $(y_1,y_2)\in I_1$, using that for $i=1,2$, in the interval $(0,y_i)$ it holds $|g(x,y_i,s)|\leq cy_i^\alpha$ and thanks to the inequalities $(y_2-y_1)^\alpha\geq\delta^\alpha y_2^\alpha\geq\delta^\alpha y_i^\alpha$, then
\begin{eqnarray*}
\frac{|\mathcal G_\varepsilon(x,y_1)-\mathcal G_\varepsilon(x,y_2)|}{(y_2-y_1)^\alpha}&\leq&\frac{1}{(y_2-y_1)^\alpha}\sum_{i=1}^2|\mathcal G_\varepsilon(x,y_i)|\\
&\leq&\frac{c}{\delta^\alpha}\sum_{i=1}^2\frac{y_i^\alpha\int_0^{y_i}\rho_\varepsilon^{-a}(s)}{y_i^\alpha\int_0^{y_i}\rho_\varepsilon^{-a}(s)}=\frac{2c}{\delta^\alpha}.
\end{eqnarray*}
If we consider $(y_1,y_2)\in I_2$, then
\begin{eqnarray*}
\frac{|\mathcal G_\varepsilon(x,y_1)-\mathcal G_\varepsilon(x,y_2)|}{(y_2-y_1)^\alpha}&\leq&\frac{1}{(y_2-y_1)^\alpha}\frac{\int_{0}^{y_2}\rho_\varepsilon^{-a}(s)|g(x,y_2,s)-g(x,y_1,s)|}{\int_0^{y_2}\rho_\varepsilon^{-a}(s)}\\
&&+\frac{1}{(y_2-y_1)^\alpha}\frac{\int_{y_1}^{y_2}\rho_\varepsilon^{-a}(s)|g(x,y_1,s)|}{\int_0^{y_2}\rho_\varepsilon^{-a}(s)}\\
&&+\frac{1}{(y_2-y_1)^\alpha}\frac{\left(\int_{0}^{y_1}\rho_\varepsilon^{-a}(s)|g(x,y_1,s)|\right)\left(\int_{y_1}^{y_2}\rho_\varepsilon^{-a}(s)\right)}{\left(\int_{0}^{y_1}\rho_\varepsilon^{-a}(s)\right)\left(\int_{0}^{y_2}\rho_\varepsilon^{-a}(s)\right)}\\
&=&J_1+J_2+J_3.
\end{eqnarray*}
Hence, $J_1$ can be bounded using the fact that $|g(x,y_2,s)-g(x,y_1,s)|\leq c(y_2-y_1)^\alpha$. Working on $J_2$, there exists $y_1\leq\xi\leq y_2$ such that
\begin{eqnarray*}
J_2&\leq& c(y_2-y_1)^{1-\alpha}\frac{\rho_\varepsilon^{-a}(\xi)y_1^\alpha}{\int_0^{y_2}\rho_\varepsilon^{-a}(s)}\\
&\leq&c\left(\frac{y_2-y_1}{y_2}\right)^{1-\alpha}\frac{y_2\rho_\varepsilon^{-a}(\xi)}{\int_0^{y_2}\rho_\varepsilon^{-a}(s)}\\
&\leq&c\delta^{1-\alpha}\max\{1,(1-\delta)^{-a}\}\frac{y_2\rho_\varepsilon^{-a}(y_2)}{\int_0^{y_2}\rho_\varepsilon^{-a}(s)}
\end{eqnarray*}
using the fact that $y_2-y_1<\delta y_2$, the inequalities
\begin{equation*}
1-\delta<\frac{y_1}{y_2}\leq\frac{\xi}{y_2}\leq 1,
\end{equation*}
and the fact that $\rho_\varepsilon^{-a}(\xi)\leq \max\{1,(1-\delta)^{-a}\}\rho_\varepsilon^{-a}(y_2)$ (easy to check). Eventually, recalling $y_2/\varepsilon=t\in[0,+\infty)$, we have already remarked that the function defined in \eqref{psi} is bounded uniformly in $\varepsilon$
$$\frac{y_2\rho_\varepsilon^{-a}(y_2)}{\int_0^{y_2}\rho_\varepsilon^{-a}(s)}=\frac{t(1+t^2)^{-a/2}}{\int_0^{t}(1+s^2)^{-a/2}}=\psi(t)\leq\max\{1,1-a\}.$$
With analogous computations we can bound also $J_3$.
\endproof

\begin{Proposition}\label{c0alphaG}
Let $a\in(-\infty,1)$, $\varepsilon\geq0$, $\alpha\in(0,1)$ and let $\gamma\in C^{0,\alpha}(B_1^+)$. Then the family of functions
\begin{equation*}
\mathcal G_\varepsilon(x,y)=\frac{\int_0^y\rho_\varepsilon^{-a}(s)\left(\gamma(x,s)-\gamma(x,0)\right)\mathrm{d}s}{\int_0^y\rho_\varepsilon^{-a}(s)\mathrm{d}s}
\end{equation*}
is uniformly bounded in $C^{0,\alpha}(B_1^+)$ by a constant which does not depend on $\varepsilon$.
\end{Proposition}
\proof
Just notice that, since $g(x,y,s):=\gamma(x,s)-\gamma(x,0)$ for $s\leq y$,  $g$ satisfies conditions of the previous Lemma \ref{a2}. Indeed is $\alpha$-H\"older continuous in $(x,y)$ uniformly in $s\leq y$ and
$$|g(x,y,s)|=|\gamma(x,s)-\gamma(x,0)|\leq c|s|^\alpha\leq c|y|^\alpha.$$
\endproof

\begin{Proposition}\label{c1alphaG}
Let $a\in(-\infty,1)$, $\varepsilon\geq0$, $\alpha\in(0,1)$ and let $\gamma\in C^{1,\alpha}(B_1^+)$ with $\partial_y\gamma(x,0)\in C^{1,\alpha}(B_1^+)$. Consider the family of functions
\begin{equation*}
\mathcal G_\varepsilon(x,y)=\frac{\int_0^y\rho_\varepsilon^{-a}(s)\left(\gamma(x,s)-\gamma(x,0)\right)\mathrm{d}s}{\int_0^y\rho_\varepsilon^{-a}(s)\mathrm{d}s}\;.
\end{equation*}
Then there exists $c>0$ such that, for every $\varepsilon\in[0,\varepsilon_0]$,
$\mathcal G_\varepsilon$ is uniformly bounded in $C^{1,\alpha}(B_1\cap\{y\geq \sqrt{\varepsilon}\})$ by $c$.

\end{Proposition}
\proof
One can rewrite our function as
\begin{eqnarray*}
\mathcal G_\varepsilon(x,y)&=&\frac{\int_0^y\rho_\varepsilon^{-a}(s)\left(\gamma(x,s)-\gamma(x,0)-\partial_y\gamma(x,0)s\right)\mathrm{d}s}{\int_0^y\rho_\varepsilon^{-a}(s)\mathrm{d}s}+\partial_y\gamma(x,0)\frac{\int_0^y\rho_\varepsilon^{-a}(s)s\, \mathrm{d}s}{\int_0^y\rho_\varepsilon^{-a}(s)\mathrm{d}s}\\
&=& \frac{\int_0^y\rho_\varepsilon^{-a}(s)\left(\int_0^s(\partial_y\gamma(x,\tau)-\partial_y\gamma(x,0))\mathrm{d}\tau\right)\mathrm{d}s}{\int_0^y\rho_\varepsilon^{-a}(s)\mathrm{d}s}+\partial_y\gamma(x,0)\frac{\int_0^y\rho_\varepsilon^{-a}(s)s\, \mathrm{d}s}{\int_0^y\rho_\varepsilon^{-a}(s)\mathrm{d}s}\\
&=& \mathcal{H}_\varepsilon(x,y)+\partial_y\gamma(x,0)\frac{\int_0^y\rho_\varepsilon^{-a}(s)s\, \mathrm{d}s}{\int_0^y\rho_\varepsilon^{-a}(s)\mathrm{d}s}.
\end{eqnarray*}
First we show that the second term has the desired property uniformly in $\varepsilon$. At first we remark that $\partial_y\gamma(x,0)\in C^{1,\alpha}(B_1^+)$. Now consider  that the family of functions
\begin{equation*}
\xi_\varepsilon^a(y):=\frac{\int_0^y\rho_\varepsilon^{-a}(s)s\, \mathrm{d}s}{\int_0^y\rho_\varepsilon^{-a}(s)\mathrm{d}s}
\end{equation*}
is uniformly bounded in $L^\infty(B_1^+)$. In fact, denoting $t=y/\varepsilon$,
\begin{equation*}
\xi_\varepsilon^a(y)=\varepsilon\xi_1^a(t)=\varepsilon\frac{\int_0^t(1+s^2)^{-a/2}s\, \mathrm{d}s}{\int_0^t(1+s^2)^{-a/2}\mathrm{d}s}=y \, \frac{\xi_1^a(t)}{t},
\end{equation*}
is bounded in $B_1^+$ (uniformly with respect to $\varepsilon\geq0$). In fact, the first factor $y$ is obviously bounded in $[0,1]$ and the second one is bounded for $t\in[0,+\infty)$. Now, let us consider the derivative in $y$,
\begin{equation*}
\partial_y\xi_\varepsilon^a(y)=(\xi_1^a)'(t)=\psi_1^a(t)\left(1-\frac{\int_0^t(1+s^2)^{-a/2}s\, \mathrm{d}s}{t\int_0^t(1+s^2)^{-a/2}\mathrm{d}s}\right)\;.
\end{equation*}

We claim that $\partial_y\xi_\varepsilon^a$ enjoys the property stated in \eqref{eq:property}, being the product of two functions, both bounded, 
having a finite limit as $t\to +\infty$ and  derivatives vanishing as $1/t^2$.
.

Eventually we consider $\mathcal H_\varepsilon$. Computing the gradient $\nabla_x\mathcal H_\varepsilon$, we obtain
\begin{equation*}
\nabla_x\mathcal H_\varepsilon(x,y)=\frac{\int_0^y\rho_\varepsilon^{-a}(s)\left(\tilde\gamma(x,s)-\tilde\gamma(x,0)\right)\mathrm{d}s}{\int_0^y\rho_\varepsilon^{-a}(s)\mathrm{d}s}
\end{equation*}
where
$$\tilde\gamma(x,s)=\nabla_x\gamma(x,s)-\nabla_x\partial_y\gamma(x,0)s\in C^{0,\alpha}(B_1^+),$$
and satisfies the assumptions in Proposition \ref{c0alphaG}.\\\\
It remains to consider the partial derivative in $y$ of $\mathcal H_\varepsilon$; that is,
\begin{eqnarray*}
\partial_y\mathcal H_\varepsilon(x,y)&=&\frac{y\rho_\varepsilon^{-a}(y)}{\int_0^y\rho_\varepsilon^{-a}(s)\mathrm{d}s}\cdot\dfrac{\int_0^y\rho_\varepsilon^{-a}(s)\left(\frac{1}{y}\int_s^y(\partial_y\gamma(x,\tau)-\partial_y\gamma(x,0))\mathrm{d}\tau\right)\mathrm{d}s}{\int_0^y\rho_\varepsilon^{-a}(s)\mathrm{d}s}\\
&=&\psi_\varepsilon^a(y)\cdot\mathcal I_\varepsilon(y).
\end{eqnarray*}
By Remark \ref{psia}, the family of functions $\psi_\varepsilon^a$ enjoy the desired propery \eqref{eq:property}. Now we wish to conclude that $\mathcal I_\varepsilon$ is uniformly bounded in $C^{0,\alpha}(B_1^+)$. To this aim, it is enough to prove that the function
$$g(x,y,s)=\frac{1}{y}\int_s^y(\partial_y\gamma(x,\tau)-\partial_y\gamma(x,0))\mathrm{d}\tau$$
satisfies conditions in Lemma \ref{a2}. Using the H\"older continuity of $\partial_y\gamma$, obviously
$$|g(x,y,s)|\leq\frac{1}{y}\int_s^y|\partial_y\gamma(x,\tau)-\partial_y\gamma(x,0)|\mathrm{d}\tau\leq\frac{c|y|^\alpha(y-s)}{y}\leq c|y|^\alpha.$$
The H\"older continuity of $g$ in the $x$-variable is trivial. Nevertheless, following the reasonings in the proof of Lemma \ref{a2}, fixed $0<\delta<1$, let us consider the following two sets
$$I_1=\{(y_1,y_2)\ : \ 0\leq y_1\leq y_2<1, \ y_2-y_1\geq\delta y_2\}$$
and
$$I_2=\{(y_1,y_2)\ : \ 0\leq y_1\leq y_2<1, \ y_2-y_1<\delta y_2\}.$$
If we consider $(y_1,y_2)\in I_1$, using that for $i=1,2$ it holds $|g(x,y_i,s)|\leq cy_i^\alpha$ and thanks to the inequalities $(y_2-y_1)^\alpha\geq\delta^\alpha y_2^\alpha\geq\delta^\alpha y_i^\alpha$, then
\begin{eqnarray*}
\frac{|g(x,y_1,s)-g(x,y_2,s)|}{(y_2-y_1)^\alpha}&\leq&\frac{1}{(y_2-y_1)^\alpha}\sum_{i=1}^2|g(x,y_i,s)|\\
&\leq&\frac{c}{\delta^\alpha}\sum_{i=1}^2\frac{y_i^\alpha}{y_i^\alpha}=\frac{2c}{\delta^\alpha}.
\end{eqnarray*}
If we consider $(y_1,y_2)\in I_2$, then, using the fact that $y_2-y_1<\delta y_2$
\begin{eqnarray*}
\frac{|g(x,y_1,s)-g(x,y_2,s)|}{(y_2-y_1)^\alpha}&\leq&\frac{1}{(y_2-y_1)^{\alpha}y_2}\int_{y_1}^{y_2}|\partial_y\gamma(x,\tau)-\partial_y\gamma(x,0)|\mathrm{d}\tau\\
&&+\frac{1}{(y_2-y_1)^{\alpha}}\left|\frac{1}{y_2}-\frac{1}{y_1}\right|\int_{s}^{y_1}|\partial_y\gamma(x,\tau)-\partial_y\gamma(x,0)|\mathrm{d}\tau\\
&\leq&c\delta^{1-\alpha}+c\frac{\delta^{1-\alpha}}{1-\delta}.
\end{eqnarray*}
\endproof

\section{Quadratic forms, stability and isometries}
In this appendix we are going to prove some useful inequalities, needed when working in weighted Sobolev spaces, specially whenever the weight does not belong to the $A_2$ class. These results will be the key of the validity of Liouville type theorems in Section \ref{sec:liouville}.
\subsection{Hardy type inequalities}
At first, we deal with  the validity of Hardy (trace) type inequalities and their spectral stability. These results will be the key tools in order to establish a class of Liouville theorems contained in this section. Let $\mathbb{R}^{n+1}_+=\mathbb{R}^{n+1}\cap\{y>0\}$, $B_1^+=B_1\cap\{y>0\}$ and $S^n_+=S^{n}\cap\{y>0\}$. We define the space $\tilde H^{1}(B_1^+)$ as the closure of $C^\infty_c(\overline B_1^+\setminus\Sigma)$ with respect to the norm
$$\left(\int_{B_1^+}|\nabla v|^2\right)^{1/2}.$$
Then, we remark that the following trace Poincar\'e inequality holds
\begin{equation}\label{tracepoinc}
c\int_{S^{n}_+}v^2\leq\int_{B_1^+}|\nabla v|^2.
\end{equation}
We first state the following Hardy inequality.
\begin{Lemma}[Hardy inequality]\label{HARDY}
Let $v\in \tilde H^{1}(B_1^+)$. Then
\begin{equation}\label{Hardy}
\frac{1}{4}\int_{B_1^+}\frac{v^2}{y^2}\leq\int_{B_1^+}|\nabla v|^2.
\end{equation}
\end{Lemma}
\proof
The proof is an easy exercise based on the well known Hardy inequality on the half space
$$\frac{1}{4}\int_{\R^{n+1}_+}\frac{v^2}{y^2}\leq\int_{\R^{n+1}_+}|\nabla v|^2,$$
and using the Kelvin transform.
\endproof

Next, we will need a boundary version of the Hardy inequality

\begin{Lemma}[Boundary Hardy inequality]\label{BoundaryHardy}
There exists $c_0>0$ such that, for every $v\in\tilde H^1(B_1^+)$, there holds
\begin{equation}\label{tracehard}
c_0\int_{S^n_+}\frac{v^2}{y}\leq\int_{B_1^+}|\nabla v|^2.
\end{equation}
\end{Lemma}
\proof
By taking the harmonic replacement of $v$ on $B_1^+$, we may assume without loss of generality that
\( \Delta v=0\) in $B_1^+$. Now we consider the following inversion (stereographic projection) $\Phi:B_1^+\subset\mathbb{R}^{n+1}\to\mathbb{R}^{n+1}$ such that
$$\Phi:z=(x,y)=(x_1,...,x_{n},y)\mapsto \tilde z=(\tilde x,\tilde y)=(\tilde x_1,...,\tilde x_{n},\tilde y),$$
with
$$\Phi(z)=\frac{z+e_1}{|z+e_1|^2}-\frac{e_1}{2}\qquad\mathrm{and}\qquad \Phi^{-1}(\tilde z)=\frac{\tilde z+\frac{e_1}{2}}{|z+\frac{e_1}{2}|^2}-e_1.$$
This map is conformal and such that $\Phi(B_1^+)=\{\tilde x_1>0\}\cap\{\tilde y>0\}$ and $\Phi(S^n_+)=\{\tilde x_1=0\}\cap\{\tilde y>0\}$. Hence, the Kelvin transform
\begin{equation*}
w(\tilde z)=Kv(\tilde z):=\frac{1}{|\tilde z+\frac{e_1}{2}|^{n-1}}v(\Phi^{-1}(\tilde z))
\end{equation*}
is harmonic in $\{\tilde x_1>0\}\cap\{\tilde y>0\}$ and such that
$$\int_{B_1^+}|\nabla v|^2\mathrm{d}z=\int_{\{\tilde x_1>0\}\cap\{\tilde y>0\}}|\nabla w|^2\mathrm{d}\tilde z.$$
Using a fractional Hardy inequality (see \cite{BogDyd}) on the $n$-dimensional half space $\{\tilde x_1=0\}\cap\{\tilde y>0\}$, up to extending the function $w=0$ in $\{\tilde x_1=0\}\cap\{\tilde y<0\}$, we have
\begin{eqnarray*}
\int_{\{\tilde x_1>0\}\cap\{\tilde y>0\}}|\nabla w|^2\mathrm{d}\tilde z &\geq& c\iint_{(\{\tilde x_1=0\}\cap\{\tilde y>0\})^2}\frac{|w(\tilde \zeta_1)-w(\tilde \zeta_2)|^2}{|\tilde \zeta_1-\tilde \zeta_2|^{n+1}}\mathrm{d}\tilde \zeta_1\mathrm{d}\tilde \zeta_2\nonumber\\
&\geq& c\int_{\{\tilde x_1=0\}\cap\{\tilde y>0\}}\frac{w^2(\tilde z)}{\tilde y}\mathrm{d}\tilde z.
\end{eqnarray*}
Finally we compute
\begin{eqnarray*}
&&\int_{S^n_+}\frac{v^2(z)}{y}\mathrm{d}\sigma(z) \\
&=&\int_{\{\tilde x_1=0\}\cap\{\tilde y>0\}}\frac{w^2(\tilde z)}{\tilde y}\left|\tilde z+\frac{e_1}{2}\right|^{2(n-1)+2} \cdot|\Phi^{-1}_{\tilde x_2}(\tilde z)\wedge\Phi^{-1}_{\tilde x_3}(\tilde z)\wedge...\wedge\Phi^{-1}_{\tilde x_{n}}(\tilde z)\wedge\Phi^{-1}_{\tilde y}(\tilde z)|\mathrm{d}\tilde z\nonumber\\
&\leq&\int_{\{\tilde x_1=0\}\cap\{\tilde y>0\}}\frac{w^2(\tilde z)}{\tilde y}\mathrm{d}\tilde z.
\end{eqnarray*}
\endproof

\subsection{A stability result}
\begin{Lemma}\label{quad}
Let $\{Q_k\}_{k\in\mathbb{N}}$ be a family of quadratic forms $Q_k:\tilde H^1(B_1^+)\to[0,+\infty)$ defined by
\begin{equation*}
Q_k(v)=\int_{B_1^+}|\nabla v|^2+\int_{B_1^+}V_kv^2+\int_{S^n_+}W_kv^2.
\end{equation*}
Assume that the family $\{Q_k\}$ satisfies the following conditions:
\begin{itemize}
\item[i)] $|W_k|\leq c$ on $S^n_+$ and $|V_k|\leq\frac{c}{y^2}$ in $B_1^+$ uniformly on $k\in\mathbb{N}$;
\item[ii)] there exists a constant $C>0$ which does not depend on $k\in\mathbb{N}$ such that for any $v\in\tilde H^1(B_1^+)$
\begin{equation}\label{unifequiv}
\frac{1}{C}\|v\|_{\tilde H^1(B_1^+)}^2\leq Q_k(v)\leq C\|v\|_{\tilde H^1(B_1^+)}^2;
\end{equation}
\item[iii)] $V_k\to V$ in $B_1^+$ and $W_k\to W$ on $S^n_+$ pointwisely as $k\to+\infty$, where
$$Q(v)=\int_{B_1^+}|\nabla v|^2+\int_{B_1^+}Vv^2+\int_{S^n_+}Wv^2,$$
with $Q:\tilde H^1(B_1^+)\to[0,+\infty)$ satisfying $|W|\leq c$ on $S^n_+$, $|V|\leq\frac{c}{y^2}$ in $B_1^+$ and 
$$\frac{1}{C}\|v\|_{\tilde H^1(B_1^+)}^2\leq Q(v)\leq C\|v\|_{\tilde H^1(B_1^+)}^2.$$
\end{itemize}
Let
$$\lambda_k=\min_{v\in\tilde H^1(B_1^+)\setminus\{0\}}\frac{Q_k(v)}{\int_{S^n_+}v^2},\qquad\lambda=\min_{v\in\tilde H^1(B_1^+)\setminus\{0\}}\frac{Q(v)}{\int_{S^n_+}v^2}.$$
Then, $\lambda_k\to\lambda$.
\end{Lemma}
\proof
Let $\{v_k\}\subset \tilde H^1(B_1^+)\setminus\{0\}$ be a sequence of minimizers for $\lambda_k$; that is, such that
$$\lambda_k=Q_{k}(v_k)=\int_{B_1^+}|\nabla v_k|^2+\int_{B_1^+}V_kv_k^2+\int_{S^n_+}W_kv_k^2,$$
and $\int_{S^n_+}v^2_k=1$. Since by compact embedding $\tilde H^1(B_1^+)\hookrightarrow L^2(S^n_+)$ the minimum
$$\min_{v\in\tilde H^1(B_1^+)\setminus\{0\}}\frac{\|v\|^2_{\tilde H^1(B_1^+)}}{\int_{S^n_+}v^2}=\frac{\|u\|^2_{\tilde H^1(B_1^+)}}{\int_{S^n_+}u^2}=\nu>0$$
is achieved by $u\in\tilde H^1(B_1^+)\setminus\{0\}$ and it is strictly positive by the trace Poincar\'e inequality, then there exists a positive constant $C$ independent from $k$ such that
$$\frac{\nu}{C}\leq\lambda_k\leq C\nu.$$
Moreover, we have that
$$\frac{1}{C}\|v_k\|^2_{\tilde H^1(B_1^+)}\leq\lambda_k\leq C\nu$$
and so there exists $\overline v\in\tilde H^1(B_1^+)$ such that $v_k\rightharpoonup\overline v$ in $\tilde H^1(B_1^+)$ and, up to passing to a subsequence, $v_k\to\overline v$ in $L^2(S^n_+)$. Moreover, the limit is non trivial by the condition $\int_{S^n_+}\overline v^2=1$.\\\\
We want to prove that the convergence is strong in $\tilde H^1(B_1^+)$. Testing the eigenvalue equation solved by $v_k$ with $v_k-\overline v$, we have
$$\int_{B_1^+}\nabla v_k\cdot\nabla(v_k-\overline v)+\int_{B_1^+}V_kv_k(v_k-\overline v)+\int_{S^n_+}W_kv_k(v_k-\overline v)=\lambda_k\int_{S^n_+}v_k(v_k-\overline v).$$
Using the fact that $|W_k|, |\lambda_k|\leq c$ uniformly in $k$, the strong convergence and the normalization in $L^2(S^n_+)$, by the H\"older inequality the terms over the half sphere $S^n_+$ go to 0 in the limit. So
\begin{equation}\label{testing11}
\int_{B_1^+}\nabla v_k\cdot\nabla(v_k-\overline v)+\int_{B_1^+}V_kv_k(v_k-\overline v)\to0.
\end{equation}
Hence,
\begin{eqnarray}\label{strong1}
Q_{k}(v_k-\overline v)&=&\int_{B_1^+}|\nabla (v_k-\overline v)|^2+\int_{B_1^+}V_k(v_k-\overline v)^2+\int_{S^n_+}W_k(v_k-\overline v)^2\nonumber\\
&=&\int_{B_1^+}\nabla v_k\cdot\nabla(v_k-\overline v)+\int_{B_1^+}V_kv_k(v_k-\overline v)-\int_{B_1^+}\nabla \overline v\cdot\nabla(v_k-\overline v)\nonumber\\
&&-\int_{B_1^+}V\overline v(v_k-\overline v)+\int_{B_1^+}(V-V_k)\overline v(v_k-\overline v)+\int_{S^n_+}W_k(v_k-\overline v)^2\to0.
\end{eqnarray}
In fact, the sum of the first two terms goes to 0 by \eqref{testing11}, the sum of the second two by weak convergence in $\tilde H^1(B_1^+)$. The third term is such that
\begin{eqnarray*}
\int_{B_1^+}(V-V_k)\overline v(v_k-\overline v)&\leq&\left(\int_{B_1^+}(V-V_k)\overline v^2\right)^{1/2}\left(\int_{B_1^+}(V-V_k)(v_k-\overline v)^2\right)^{1/2}\\
&\leq& c\left(\int_{B_1^+}(V-V_k)\overline v^2\right)^{1/2}\to0.
\end{eqnarray*}
We used that $V_k\to V$, the fact that $|V_k-V|\leq\frac{c}{y^2}$ and the Hardy inequality to ensure the dominated convergence theorem. Eventually the last term in \eqref{strong1} goes to 0 by the strong convergence in $L^2(S^n_+)$. Hence we obtain the strong convergence by \eqref{unifequiv}.\\\\
It is easy to see that $Q_{k}(v_k)\to Q(\overline v)$. This is enough to conclude because if we consider $\tilde v$ the normalized in $L^2(S^n_+)$ minimizer of $\lambda$, since it is competitor for the minimization of any $Q_{k}$, then
$$\lambda_k=Q_{k}(v_k)\leq Q_{k}(\tilde v),$$
and since $Q_{k}(v_k)\to Q(\overline v)$ and $Q_{k}(\tilde v)\to Q(\tilde v)$, then by $Q(\overline v)\leq Q(\tilde v)$, and by the minimality of $\overline v$, we finally obtain that $\overline v=\tilde v$ with $\lambda_k\to\lambda$.
\endproof

\subsection{Quadratic forms for the odd case}
Let $a\in(-\infty,1)$, $\varepsilon\geq0$ and consider a function $u\in C^\infty_c(\overline B_1^+\setminus\Sigma)$ and define $v=(\rho_\varepsilon^a)^{1/2}u\in C^\infty_c(\overline B_1^+\setminus\Sigma)$. Let us define the quadratic form
\begin{equation}\label{Qrho}
\int_{B_1^+}\rho_\varepsilon^au^2=Q_{\rho_\varepsilon^a}(v)=\int_{B_1^+}|\nabla v|^2+\int_{B_1^+}V_{\rho_\varepsilon^a}v^2+\int_{S^n_+}W_{\rho_\varepsilon^a}v^2,
\end{equation}
where
$$V_{\rho_\varepsilon^a}(y)=\frac{(\rho_\varepsilon^a)''}{2\rho_\varepsilon^a}-\left(\frac{(\rho_\varepsilon^a)'}{2\rho_\varepsilon^a}\right)^2=\frac{a[(a-2)y^2+2\varepsilon^2]}{4(\varepsilon^2+y^2)^2}$$
and
$$W_{\rho_\varepsilon^a}(y)=-\frac{(\rho_\varepsilon^a)'y}{2\rho_\varepsilon^a}=-\frac{ay^2}{2(\varepsilon^2+y^2)}.$$
Let 
\begin{equation*}
Q_{a}(v)=\int_{B_1^+}|\nabla v|^2+\int_{B_1^+}V_{a}v^2+\int_{S^n_+}W_{a}v^2,
\end{equation*}
with $V_a(y)=\frac{a(a-2)}{4y^2}=V_{\rho_0^a}(y)$ and $W_{a}(y)=-\frac{a}{2}=W_{\rho_0^a}(y)$. Eventually consider a sequence $\varepsilon_k\to0$ as $k\to+\infty$ and define $\rho_k=\rho_{\varepsilon_k}^a$. Let us recall $Q_k=Q_{\rho_k}$ and $Q=Q_a$.
\begin{Lemma}\label{A1}
Under the previous hypothesis, the family $\{Q_{k}\}=\{Q_{\rho_{\varepsilon_k}}\}$ defined in \eqref{Qrho} and its limit $Q$ satisfy the conditions in Lemma \ref{quad}.
\end{Lemma}
\proof
Condition $i)$ is trivially satisfied. Moreover, combining $i)$, the trace Poincar\'e and the Hardy inequalities, we easily obtain the upper bound in $ii)$ for any $k\in\mathbb{N}$ with a constant independent on $\varepsilon_k$; that is,
$$Q_k(v)\leq c\|v\|_{\tilde H^1(B_1^+)}^2.$$
Let us consider $Q=Q_a$ and let us define $u=y^{-a/2}v\in C^\infty_c(\overline B_1^+\setminus\Sigma)$.
\begin{eqnarray}\label{Qa}
Q_a(v)&=&\int_{B_1^+}|\nabla v|^2+\left(\frac{a^2}{4}-\frac{a}{2}\right)\int_{B_1^+}\frac{v^2}{y^2}-\frac{a}{2}\int_{S^n_+}v^2\\\nonumber
&=&\int_{B_1^+}|\nabla v|^2+\left(\frac{a^2}{4}-\frac{a}{2}\right)\int_{B_1^+}\frac{v^2}{y^2}-\frac{a}{2}\int_{B_1^+}\mathrm{div}\left(\frac{v^2}{y}\vec{e_n}\right)=\int_{B_1^+}y^a|\nabla u|^2.
\end{eqnarray}
First of all we notice that if $a\leq 0$ the lower bound follows trivially. So we can suppose that $a\in(0,1)$. Since for $a\neq1$,  $(\frac{a^2}{4}-\frac{a}{2})>-\frac{1}{4}$, hence by the Hardy inequality in \eqref{Hardy}, the quantity
$$G_a(v)=\int_{B_1^+}|\nabla v|^2+\left(\frac{a^2}{4}-\frac{a}{2}\right)\int_{B_1^+}\frac{v^2}{y^2}$$
defines an equivalent norm in $\tilde H^1(B_1^+)$. Hence by the compact embedding $\tilde H^1(B_1^+)\hookrightarrow L^2(S^n_+)$ we have that the minimum in
$$\xi(a)=\min_{v\in\tilde H^1(B_1^+)\setminus\{0\}}\frac{G_a(v)}{\int_{S^n_+}v^2}$$
is achieved. In fact, considering a minimizing sequence, we can take it such that $\int_{S^n_+}v_k^2=1$ and also such that $v_k\in C^\infty_c(\overline B_1^+\setminus\Sigma)$. So it is uniformly bounded in $\tilde H^1(B_1^+)$ and $v_k\rightharpoonup\overline v\in\tilde H^1(B_1^+)$ with $G_a(v_k)\to\xi(a)$. Moreover the convergence is strong in $L^2(S^n_+)$ by compact embedding. Since $\int_{S^n_+}v_k^2=1$, we also obtain convergence of the $\tilde H^1_0$-norms of the $v_k$ to that of the limit, yielding strong convergence in $\tilde H^1(B_1^+)$. In fact, by the lower semicontinuity of the norm
$$\xi(a)\leq\frac{G_a(\overline v)}{\int_{S^n_+}\overline v^2}\leq\liminf_{k\to+\infty}\frac{G_a(v_k)}{\int_{S^n_+}v_k^2}=\xi(a).$$
Obviously by the condition $\int_{S^n_+}\overline v^2=1$ the limit $\overline v$ is not trivial. This proves that $\overline v$ achieves the minimum. Moreover, defining
\begin{equation}\label{lama}
\lambda(a)=\min_{v\in\tilde H^1(B_1^+)\setminus\{0\}}\frac{Q_a(v)}{\int_{S^n_+}v^2}=\xi(a)-\frac{a}{2}\geq0,
\end{equation}
we want to prove that actually $\lambda(a)>0$. First of all, such a minimum is nonnegative since the minimizing sequence can be taken in $C^\infty_c(\overline B_1^+\setminus\Sigma)$ and so the equalities in \eqref{Qa} give this condition. By contradiction let $\lambda(a)=0$. Hence the minimizing sequence is such that $Q_a(v_k)\to0$. Defining $u_k=y^{-a/2}v_k$, one has $\int_{B_1^+}y^a|\nabla u_k|^2\to0$. Moreover, the strong convergence in $\tilde H^1(B_1^+)$ gives the almost everywhere convergence of $\nabla v_k\to\nabla\overline v$ which of course implies that $\nabla u_k\to \nabla(y^{-a/2}\overline v)$ almost everywhere in $B_1^+$. Hence, since $\nabla(y^{-a/2}\overline v)=0$ almost everywhere, $\overline v=cy^{a/2}$, but $\nabla\overline v$ does not belong to $L^2(B_1^+)$. This is a contradiction. So $\lambda(a)>0$. So we have the inequality
$$Q_a(v)\geq\lambda(a)\int_{S^n_+}v^2,$$
which says that
$$Q_a(v)\geq\frac{\lambda(a)}{\frac{a}{2}+\lambda(a)}\left(\int_{B_1^+}|\nabla v|^2+\left(\frac{a^2}{4}-\frac{a}{2}\right)\int_{B_1^+}\frac{v^2}{y^2}\right),$$
and by the equivalence of the norms we obtain the result for a constant which depends on $a$ and $\lambda(a)$. Eventually, we have proved that also $Q_a$ is an equivalent norm on $\tilde H^1(B_1^+)$.\\\\
In order to prove the lower bound for $Q_k$ which is uniform in $k$, it is enough to remark that if $a\geq0$, then $Q_k\geq Q_a$. If $a<0$, then one can check that
$$Q_k(v)\geq \int_{B_1^+}|\nabla v|^2-\int_{B_1^+}\frac{a}{4(a-4)}\frac{v^2}{y^2},$$
with $\frac{a}{4(a-4)}<\frac{1}{4}$ and hence by the Hardy inequality in \eqref{Hardy} we have also in this case an equivalent norm.
\endproof

Let us recall the definition of $\tilde H^1(B_1^+,\rho_\varepsilon^a(y)\mathrm{d}z)$ as the closure of $C^\infty_c(\overline B_1^+\setminus\Sigma)$ with respect to the norm
$$\int_{B_1^+}\rho_\varepsilon^a|\nabla u|^2.$$
\begin{Lemma}\label{A2}
Let $a\in(-\infty,1)$, $\varepsilon\geq0$ and $u\in\tilde H^1(B_1^+,\rho_\varepsilon^a(y)\mathrm{d}z)$. Then the following inequalities hold true for a positive constant $c$ independent of $\varepsilon\in[0,1]$
\begin{equation}\label{poin}
c\int_{B_1^+}\rho_\varepsilon^au^2\leq\int_{B_1^+}\rho_\varepsilon^a|\nabla u|^2,
\end{equation}
\begin{equation}\label{tracepoin}
c\int_{S^n_+}\rho^a_{\varepsilon}u^2\leq\int_{B_1^+}\rho_\varepsilon^a|\nabla u|^2,
\end{equation}
\begin{equation}\label{hard}
c\int_{B_1^+}\frac{\rho_\varepsilon^a}{y^2}u^2\leq\int_{B_1^+}\rho_\varepsilon^a|\nabla u|^2,
\end{equation}
\begin{equation}\label{tracehard}
c\int_{S^n_+}\frac{\rho_\varepsilon^a}{y}u^2\leq\int_{B_1^+}\rho_\varepsilon^a|\nabla u|^2,
\end{equation}
\begin{equation}\label{sob}
\left(\int_{B_1^+}(\rho_\varepsilon^a)^{2^*/2}|u|^{2^*}\right)^{2/2^*}\leq c\int_{B_1^+}\rho_\varepsilon^a|\nabla u|^2,
\end{equation}
which are respectively the Poincar\'e inequality, the trace Poincar\'e inequality, the Hardy inequality, the trace Hardy inequality and a Sobolev type inequality.
\end{Lemma}
\proof
The proof is performed for functions $u\in C^\infty_c(\overline B_1^+\setminus\Sigma)$ and then extending the inqualities to $u\in\tilde H^1(B_1^+,\rho_\varepsilon^a(y)\mathrm{d}z)$ by a density argument. By Lemma \ref{A1} there exists a positive constant uniform in $\varepsilon$ such that
\begin{equation}\label{Qr}
\int_{B_1^+}\rho_\varepsilon^a|\nabla u|^2=Q_{\rho_\varepsilon^a}((\rho_\varepsilon^a)^{1/2}u)\geq c\int_{B_1^+}|\nabla((\rho_\varepsilon^a)^{1/2}u)|^2,
\end{equation}
then all the inequalities are obtained by the validity of them in $\tilde H^1(B_1^+)$.
\endproof

\subsection{Quadratic forms for the auxiliary weights}

Consider now $a\in(-\infty,1)$ and define
$$\pi_\varepsilon^a(y)=\left((1-a)\int_0^y\rho_\varepsilon^{-a}(s)\mathrm{d}s\right)^2,$$
and
\begin{equation*}
\omega_\varepsilon^a(y)=\rho_\varepsilon^a(y)\pi_\varepsilon^a(y).
\end{equation*}
We observe that this weight is super degenerate; that is, at $\Sigma$
$$\omega_\varepsilon^a(y)\sim\begin{cases}
|y|^{2-a} & \mathrm{if \ }\varepsilon=0\\
|y|^2 & \mathrm{if \ }\varepsilon>0,
\end{cases}$$
with $2-a\in(1,+\infty).$

\subsubsection{Super singular weights $(\omega_\varepsilon^a)^{-1}$}
Let us consider $u\in C^\infty_c(\overline B_1^+\setminus\Sigma)$ and define $v=(\omega_\varepsilon^a)^{-1/2}u\in C^\infty_c(\overline B_1^+\setminus\Sigma)$. Then we consider the quadratic form
\begin{equation}\label{Qomega}
\int_{B_1^+}(\omega_\varepsilon^a)^{-1}|\nabla u|^2=Q_{\omega_\varepsilon^a}(v)=\int_{B_1^+}|\nabla v|^2+\int_{B_1^+}V_{\omega_\varepsilon^a} v^2+\int_{S^n_+}W_{\omega_\varepsilon^a} v^2,
\end{equation}
with
$$V_{\omega_\varepsilon^a}=\frac{1}{4}[(\log\omega_\varepsilon^a)']^2-\frac{1}{2}(\log\omega_\varepsilon^a)'',$$
and
$$W_{\omega_\varepsilon^a}=\frac{1}{2}(\log\omega_\varepsilon^a)'y.$$
Hence
$$V_{\omega_0^a}(y)=\frac{(2-a)(4-a)}{4y^2},\qquad\mathrm{and}\qquad W_{\omega_0^a}(y)=\frac{2-a}{2}.$$
Eventually consider a sequence $\varepsilon_k\to0$ as $k\to+\infty$ and define $\omega_k=\omega_{\varepsilon_k}^a$. Let us name $Q_k=Q_{\omega_k}$ and $Q=Q_{\omega_0^a}$.
In what follows it would be useful to consider for $t>0$, the continuous function defined in \eqref{psi}; that is,
\begin{equation*}
\psi(t)=\frac{t(1+t^2)^{-a/2}}{\int_0^{t}(1+s^2)^{-a/2}\mathrm{d}s},
\end{equation*}
which is monotone nondecreasing if $a<0$ and nonincreasing if $a\in(0,1)$. Since $\psi$ has limit $1$ as $t\to0$ and limit $1-a$ as $t\to+\infty$, then
$$\sup_{t>0}\psi(t)=\max\{1,1-a\}\qquad\mathrm{and}\qquad \inf_{t>0}\psi(t)=\min\{1,1-a\}.$$
Let us finally define for any $k\in\mathbb{N}$
\begin{equation}\label{Qomega1}
\tilde Q_k(v)=Q_{k}(v)+\left(-\frac{a}{2}\right)^+\int_{S^n_+}v^2.
\end{equation}

First we need the following technical result.
\begin{Lemma}\label{Phia}
Let us define for $a\in(-\infty,1)$ and $t\in[0,+\infty)$ the function
\begin{equation}\label{phia}
\Phi_a(t)=\left[\frac{\sqrt{2}t(1+t^2)^{-a/2}}{\int_0^t(1+s^2)^{-a/2}}+\frac{at^2}{\sqrt{2}(1+t^2)}\right]^2+\frac{at^2[(2-a)t^2-2]}{4(1+t^2)^2}.
\end{equation}
Hence there exists a positive constant $c_1(a)>-\frac{1}{4}$ such that
\begin{equation}
\inf_{t>0}\Phi_a(t)=c_1(a).
\end{equation}
\end{Lemma}
\proof
{\bf Step 1: $a\in(-3,1)$.}\\\\
Whenever $0\leq a<1$, there holds
$$\min_{t>0}f_a(t)=\min_{t>0}\frac{at^2[(2-a)t^2-2]}{4(1+t^2)^2}=f_a\left(\frac{1}{\sqrt{3-a}}\right)=\frac{a}{4(a-4)}>-\frac{1}{4}.$$
Moreover, if $a<0$, 
$$\inf_{t>0}f_a(t)=\lim_{t\to+\infty}f_a(t)=\frac{a(2-a)}{4}.$$
Hence, whenever $1-\sqrt{2}<a<0$, then, the infimum remains strictly larger that $-1/4$.\\\\
Moreover, for $a<0$, then $f_a(t)\geq0$ in $\left[0,\sqrt{\frac{2}{2-a}} \ \right]$. From now on we will consider $a<0$ and $t>\sqrt{\frac{2}{2-a}} $. Now, let us compute the square in \eqref{phia}, and add $1/4$; that is
\begin{eqnarray*}
\Phi_a(t)+\frac{1}{4}&=&\frac{2t^2(1+t^2)^{-a}}{\left(\int_0^t(1+s^2)^{-a/2}\right)^2}+\frac{2at^3(1+t^2)^{-1-a/2}}{\int_0^t(1+s^2)^{-a/2}}+\frac{a^2t^4}{2(1+t^2)^2}+f_a(t)+\frac{1}{4}\\
&=&\frac{2t^3(1+t^2)^{-1-a/2}}{\left(\int_0^t(1+s^2)^{-a/2}\right)^2} \ \cdot g_a(t)+\frac{t^4(a^2+2a+1)+t^2(-2a+2)+1}{4(1+t^2)^2}\\
&=&I_a(t)+J_a(t),
\end{eqnarray*}
with
$$g_a(t)=\left(\frac{(1+t^2)^{1-a/2}}{t}+a\int_0^t(1+s^2)^{-a/2}\right).$$
It is easy to see that
$$\inf_{t>0}J_a(t)\begin{cases}>0 & \mathrm{if \ } a\neq-1\\ =0 & \mathrm{if \ }a=-1.\end{cases}$$
Nevertheless, since
$$g_a'(t)=\frac{(1+t^2)^{-a/2}}{t^2}(t^2-1),$$
then $g_a$ has its global minimum in $t=1$, and hence it is easy to see that
$$g_a(1)=2^{1-a/2}+a\int_0^1(1+s^2)^{-a/2}\geq 2^{1-a/2}+a\int_0^1(1+s)^{-a/2}=2^{1-a/2}\frac{2+a}{2-a}-\frac{2a}{2-a}>0,$$
surely if $a>-3$. Hence, when $a\in(-3,-1)\cup(-1,0)$, we have the result since $\inf_{t>0}I_a(t)\geq0$ and $\inf_{t>0}J_a(t)>0$. In the case $a=-1$ one can see that
$$\inf_{t>0}I_{-1}(t)=\min_{t>0}I_{-1}(t)>0,$$
using the explicit form
$$I_{-1}(t)=\frac{2t^3(1+t^2)^{-1-a/2}}{\frac{1}{4}\left(t\sqrt{t^2+1}+\log(\sqrt{t^2+1}+t)\right)^2}\left(\frac{(1+t^2)^{1-a/2}}{t}-\frac{1}{2}\left(t\sqrt{t^2+1}+\log(\sqrt{t^2+1}+t)\right)\right).$$
{\bf Step 2: $a\leq-3$.}\\\\
We can express 
$$\Phi_a(t)+\frac{1}{4}=\frac{t^4}{(1+t^2)^2}\left(2\left(\frac{(1+t^2)^{-\frac{a}{2}+1}}{t\int_0^t(1+s^2)^{-\frac{a}{2}}}+\frac{a}{2}\right)^2+\frac{a(2-a)}{4}-\frac{a}{2t^2}+\frac{1}{4}\frac{(1+t^2)^2}{t^4}\right).$$
Hence
$$\tilde\Phi_a(t)=\frac{(1+t^2)^2}{t^4}\left(\Phi_a(t)+\frac{1}{4}\right),$$
and $\gamma_a(t)=\tilde\Phi_a(t/\sqrt{-a})-\frac{0.001}{4}\frac{(-a+t^2)^2}{t^4}$; that is,
\begin{equation}\label{gammaa}
\gamma_a(t)=2a^2\left(\frac{(1+\frac{t^2}{-a})^{-\frac{a}{2}+1}}{t\int_0^t(1+\frac{s^2}{-a})^{-\frac{a}{2}}}-\frac{1}{2}\right)^2+\frac{a(2-a)}{4}+\frac{a^2}{2t^2}+\frac{0.999}{4}\frac{(-a+t^2)^2}{t^4}.
\end{equation}
First need to highlight some fundamental properties of the functions
\begin{equation*}
w_a(t)=\frac{(1+\frac{t^2}{-a})^{-\frac{a}{2}+1}}{t\int_0^t(1+\frac{s^2}{-a})^{-\frac{a}{2}}}.
\end{equation*}
As $a\to-\infty$ one has the pointwise convergence $w_a(t)\to v(t)$ in $(0,+\infty)$ (which is however uniform on compact subsets) with
\begin{equation*}
v(t)=\frac{e^{\frac{t^2}{2}}}{t\int_0^te^{\frac{s^2}{2}}}.
\end{equation*}
We wish to prove the following\\\\
{\bf Claim:} $w_a/v\geq 1$ in $[0,+\infty)$. At first, elementary computations show that, in a neighbourhood of $t=0$, the expansion 
$$w_a(t)=\frac{1}{t^2}+\frac{1}{2}+\frac{1}{-a}+o(1)\qquad\mathrm{and}\qquad v(t)=\frac{1}{t^2}+\frac{1}{2}+o(1),$$
holds, while in a neighbourhood of $t=+\infty$ we have
$$w_a(t)=\frac{1-a}{-a}+o(1)\qquad\mathrm{and}\qquad v(t)=1+o(1),$$
implying that $w_a/v>1$ near zero and at infinity. Thus, the claim is false if and only if there exists $t_0>0$ such that
\begin{equation}\label{syst0}
\begin{cases}
w_a(t_0)=v_a(t_0)\\
\left(\frac{w_a}{v}\right)'(t_0)\leq0.
\end{cases}
\end{equation}
Remark that, $w_a$ and $v$ solve respectively the following differential equations
\begin{equation*}
w_a'(t)=\frac{1}{t(1+\frac{t^2}{-a})}\left(\frac{1-a}{-a}t^2-1\right)w_a(t)-\frac{t}{1+\frac{t^2}{-a}}w_a^2(t)
\end{equation*}
and
\begin{equation*}
v'(t)=\frac{t^2-1}{t}v(t)-tv^2(t).
\end{equation*}
Using these equations we obtain
\begin{equation*}
\left(\frac{w_a}{v}\right)'=\frac{w_a}{v}\left(\frac{t}{-a+t^2}(2-t^2)-\frac{t}{1+\frac{t^2}{-a}}w_a+tv\right),
\end{equation*}
and \eqref{syst0} holds if and only if
\begin{equation*}
v(t_0)\leq 1-\frac{2}{t_0^2}.
\end{equation*}
Now we are going to show that, on the contrary,
\begin{equation}\label{v>z}
v(t)>z(t):=1-\frac{2}{t^2}.
\end{equation}
In $(0,\sqrt 2)$ we have $v>0$ and $z<0$. Moreove the inequality \eqref{v>z} can be checked numerically (with error estimate) on $[\sqrt 2,\sqrt 6]$, and is also valid in a neighbourhood of $t=+\infty$, by the exapnsion
$$v(t)=1-\frac{1}{t^2}+o\left(\frac{1}{t^2}\right)>z(t).$$
So, the function $v-z$ is positive near $0$ and at $+\infty$, and hence denying \eqref{v>z}  yields the existence of $t_1\geq\sqrt 6$ such that
$$\begin{cases}
v(t_1)=z(t_1)\\
(v-z)'(t_1)\leq0.
\end{cases}$$
It is easy to see that at such a point $t_1$ one has $(v-z)'(t_1)>0$ if $t_1\geq\sqrt 6$ (using the fact that $v(t_1)=z(t_1)$).\\\\
Now we can turn back to \eqref{gammaa}, obtaining by convexity that
\begin{equation}\label{gammaaa}
\gamma_a(t)\geq2a^2\left(v(t)-\frac{1}{2}\right)^2+\frac{a(2-a)}{4}+\frac{a^2}{2t^2}+\frac{0.999}{4}\frac{(-a+t^2)^2}{t^4}.
\end{equation}

In order to complete the proof, we need to prove positivity of the right hand side. To this aim, we observe that the function $v$ changes monotonicity only once on $(0,+\infty)$ and its absolute minimum value $0,77836\pm 10^{-5}$ is larger than  $1/2$. Moreover, as $v(5.1)=0.95774\pm10^{-5}$ and $v'(5.1)=0.001860\pm 10^{-5}>0$ we infer
positivity of the right hand side for  all $t\in [5.1,+\infty)$, for all $a\leq -2.96767$. The remaining values $(a,t)$ lay in the compact rectangle $[-43.3272,-2.96767]\times[1,5.1]$ and can be easily dealt numerically with error controlled minimization.

\endproof

\begin{Lemma}\label{B1}
Under the previous hypothesis, the family $\{\tilde Q_{k}\}=\{\tilde Q_{\omega_{\varepsilon_k}}\}$ defined in \eqref{Qomega1} and its limit $\tilde Q$ satisfy the conditions in Lemma \ref{quad}.
\end{Lemma}
\proof
First, we want to prove property $i)$; that is, there exists a positive constant $c>0$ uniform in $\varepsilon\to0$ such that
$$|V_{\omega_\varepsilon^a}(y)|\leq\frac{c}{y^2}\qquad\mathrm{and}\qquad |W_{\omega_\varepsilon^a}(y)|\leq c.$$
We remark that there exists a positive constant $c>0$ uniform in $\varepsilon\to0$ such that
$$|(\log\rho_\varepsilon^a)'|=\left|\frac{(\rho_\varepsilon^a)'}{\rho_\varepsilon^a}\right|\leq |a|\frac{y}{\varepsilon^2+y^2}\leq\frac{c}{y}.$$
Moreover
$$|(\log\rho_\varepsilon^a)''|\leq \left|\frac{(\rho_\varepsilon^a)'}{\rho_\varepsilon^a}\right|\cdot\left|\frac{(\rho_\varepsilon^a)''}{(\rho_\varepsilon^a)'}\right|+\left|\frac{(\rho_\varepsilon^a)'}{\rho_\varepsilon^a}\right|^2\leq \frac{c}{y^2}.$$
It remains to prove the following uniform bounds
\begin{equation*}
\left|\frac{(\pi_\varepsilon^a)'}{\pi_\varepsilon^a}\right|\leq \frac{c}{y},\qquad\mathrm{and}\qquad\left|\frac{(\pi_\varepsilon^a)''}{(\pi_\varepsilon^a)'}\right|\leq \frac{c}{y}.
\end{equation*}
Then the result follows since we are considering the logarithm of a product by linearity of the derivative.
\begin{eqnarray*}
|(\log\pi_\varepsilon^a)'|=\left|\frac{(\pi_\varepsilon^a)'}{\pi_\varepsilon^a}\right|&=&2\frac{\rho_\varepsilon^{-a}(y)}{\int_0^y\rho_\varepsilon^{-a}(s)\mathrm{d}s}\\
&=&\frac{2}{y}\frac{\frac{y}{\varepsilon}(1+\left(\frac{y}{\varepsilon}\right)^2)^{-a/2}}{\int_0^{\frac{y}{\varepsilon}}(1+s^2)^{-a/2}\mathrm{d}s}\\
&\leq&\frac{2}{y}\sup_{t>0}\frac{t(1+t^2)^{-a/2}}{\int_0^{t}(1+s^2)^{-a/2}\mathrm{d}s}\leq \frac{2\max\{1,1-a\}}{y}.
\end{eqnarray*}
Moreover,
$$\left|\frac{(\pi_\varepsilon^a)''}{(\pi_\varepsilon^a)'}\right|\leq\frac{\rho_\varepsilon^{-a}(y)}{\int_0^y\rho_\varepsilon^{-a}(s)\mathrm{d}s}+|a|\frac{y}{\varepsilon^2+y^2}\leq\frac{\max\{1,1-a\}+|a|}{y}.$$
Eventually
$$|(\log\pi_\varepsilon^a)''|\leq \left|\frac{(\pi_\varepsilon^a)'}{\pi_\varepsilon^a}\right|\cdot\left|\frac{(\pi_\varepsilon^a)''}{(\pi_\varepsilon^a)'}\right|+\left|\frac{(\pi_\varepsilon^a)'}{\pi_\varepsilon^a}\right|^2\leq \frac{c}{y^2}.$$
Obviously, point $i)$ implies the uniform upper bound in \eqref{unifequiv} by trace Poincar\'e and the Hardy inequalities. In order to prove the uniform lower bound and eventually proving $ii)$, we only have to prove that there exists a positive constant $c_1>-\frac{1}{4}$ uniform in $\varepsilon\to0$ such that
$$V_{\omega_\varepsilon^a}(y)\geq\frac{c_1}{y^2}.$$
In fact,
$$W_{\omega_\varepsilon^a}(y)+\left(-\frac{a}{2}\right)^+\geq0.$$
Let $t=y/\varepsilon>0$. Then
$$V_{\omega_\varepsilon^a}(y)=\frac{\Phi_a(t)}{y^2},$$
with $\Phi_a$ as in definition \eqref{phia}. We can conclude by applying Lemma \ref{Phia}.\\\\
Eventually we remark that also condition $iii)$ holds true.
\endproof

Let us define $\tilde H^1(B_1^+,(\omega_\varepsilon^a(y))^{-1}\mathrm{d}z)$ as the closure of $C^\infty_c(\overline B_1^+)$ with respect to the norm
$$\int_{B_1^+}(\omega_\varepsilon^a)^{-1}|\nabla u|^2.$$
\begin{Lemma}\label{A2BH}
Let $a\in(-\infty,1)$ and $u\in\tilde H^1(B_1^+,(\omega_\varepsilon^a(y))^{-1}\mathrm{d}z)$. Then 
the following inequalities hold true for a positive constant $c$ independent of $\varepsilon\in[0,1]$
\begin{equation}\label{poinBH}
c\int_{B_1^+}(\omega_\varepsilon^a)^{-1}u^2\leq\int_{B_1^+}(\omega_\varepsilon^a)^{-1}|\nabla u|^2,
\end{equation}
\begin{equation}\label{tracepoinBH}
c\int_{S^n_+}(\omega_\varepsilon^a)^{-1}u^2\leq\int_{B_1^+}(\omega_\varepsilon^a)^{-1}|\nabla u|^2,
\end{equation}
\begin{equation}\label{hardBH}
c\int_{B_1^+}\frac{(\omega_\varepsilon^a)^{-1}}{y^2}u^2\leq\int_{B_1^+}(\omega_\varepsilon^a)^{-1}|\nabla u|^2,
\end{equation}
\begin{equation}\label{tracehardBH}
c\int_{S^n_+}\frac{(\omega_\varepsilon^a)^{-1}}{y}u^2\leq\int_{B_1^+}(\omega_\varepsilon^a)^{-1}|\nabla u|^2,
\end{equation}
\begin{equation}\label{sobBH}
c\left(\int_{B_1^+}((\omega_\varepsilon^a)^{-1})^{2^*/2}|u|^{2^*}\right)^{2/2^*}\leq \int_{B_1^+}(\omega_\varepsilon^a)^{-1}|\nabla u|^2,
\end{equation}
which are respectively the Poincar\'e inequality, the trace Poincar\'e inequality, the Hardy inequality, the trace Hardy inequality and a Sobolev type inequality.
\end{Lemma}
\proof
First, we prove \eqref{tracepoinBH}. Thanks to Lemma \ref{B1} we can define for a sequence $\varepsilon_k\to0$
$$\tilde\mu_k=\min_{v\in\tilde H^1(B_1^+)\setminus\{0\}}\frac{\tilde Q_{k}(v)}{\int_{S^n_+}v^2}=\min_{v\in\tilde H^1(B_1^+)\setminus\{0\}}\frac{Q_{k}(v)}{\int_{S^n_+}v^2}+\left(-\frac{a}{2}\right)^+=\mu_k+\left(-\frac{a}{2}\right)^+,$$
and
$$\tilde\mu=\min_{v\in\tilde H^1(B_1^+)\setminus\{0\}}\frac{\tilde Q(v)}{\int_{S^n_+}v^2}=\min_{v\in\tilde H^1(B_1^+)\setminus\{0\}}\frac{Q(v)}{\int_{S^n_+}v^2}+\left(-\frac{a}{2}\right)^+=\mu+\left(-\frac{a}{2}\right)^+.$$
Actually, we are able to provide the value of $\mu$ since $u(x,y)=y^{1-(a-2)}$ is the unique function in $\tilde H^{1,a-2}(B_1^+)\setminus\{0\}$ which solves
\begin{equation*}
\begin{cases}
-L_{a-2}u=0 &\mathrm{in} \ B_1^+\\
u>0 &\mathrm{in} \ B_1^+\\
u(x,0)=0\\
\nabla u\cdot\nu=\mu u &\mathrm{in} \ S_+^{n},
\end{cases}
\end{equation*}
with $\mu=1-(a-2)=3-a$. Hence, by Lemma \ref{quad}, since $\tilde\mu_k\to\tilde\mu$, then $\mu_k\to\mu=3-a>0$ and one can find $\varepsilon_0>0$ such that for $0\leq\varepsilon_k\leq\varepsilon_{\overline k}=\varepsilon_0$ one has $\mu_k\geq\mu_{\overline k}>0$. Hence one has \eqref{tracepoinBH} with a constant $\mu_{\overline k}>0$ uniform in $0\leq\varepsilon\leq\varepsilon_0$.
For the other inequalities, the proof is done taking functions $u\in C^\infty_c(\overline B_1^+\setminus\Sigma)$ and then passing to functions $u\in\tilde H^1(B_1^+,(\omega_\varepsilon^a(y))^{-1}\mathrm{d}z)$ by density. By Lemma \ref{B1} there exists a positive constant uniform in $\varepsilon$ such that
\begin{equation}\label{QrBH}
\int_{B_1^+}(\omega_\varepsilon^a)^{-1}|\nabla u|^2+\left(-\frac{a}{2}\right)^+\int_{S^n_+}(\omega_\varepsilon^a)^{-1}u^2=\tilde Q_{\omega_\varepsilon^a}((\omega_\varepsilon^a)^{-1/2}u)\geq c\int_{B_1^+}|\nabla((\omega_\varepsilon^a)^{-1/2}u)|^2,
\end{equation}
then all the inequalities are obtained by the validity of them in $\tilde H^1(B_1^+)$ and using the trace Poincar\'e inequality \eqref{tracepoinBH}.
\endproof

\subsubsection{Super degenerate weights $\omega_\varepsilon^a$}
Let $a\in(-\infty,1)$ and let us consider $u\in C^\infty(B_1^+)$ and define $v=(\omega_\varepsilon^a)^{1/2}u\in C^\infty_c(\overline B_1^+\setminus\Sigma)$. Then we consider the quadratic form
\begin{equation}\label{Qomega1+}
\int_{B_1^+}\omega_\varepsilon^a\left(|\nabla u|^2+u^2\right)=\overline Q_{\omega_\varepsilon^a}(v)=\int_{B_1^+}\left(|\nabla v|^2+v^2\right)+\int_{B_1^+}\overline V_{\omega_\varepsilon^a} v^2+\int_{S^n_+}\overline W_{\omega_\varepsilon^a} v^2,
\end{equation}
with
$$\overline V_{\omega_\varepsilon^a}=\frac{1}{4}[(\log\omega_\varepsilon^a)']^2+\frac{1}{2}(\log\omega_\varepsilon^a)''=\frac{a}{4}\frac{(a-2)y^2+2\varepsilon^2}{(\varepsilon^2+y^2)^2},$$
and
$$\overline W_{\omega_\varepsilon^a}=-\frac{1}{2}(\log\omega_\varepsilon^a)'y.$$
We remark that $\overline V_{\omega_\varepsilon^a}=V_{\rho_\varepsilon^a}$ in \eqref{Qrho}. It is easy to check that the family of quadratic forms $\overline Q_{\omega_\varepsilon^a}$ are equivalent norms in $\tilde H^1(B_1^+)$ with constants which are uniform in $\varepsilon$; i.e. the following holds
\begin{Lemma}\label{B1+}
Under the previous hypothesis, the family $\{\overline Q_{k}\}=\{\overline Q_{\omega_{\varepsilon_k}}\}$ defined in \eqref{Qomega1+} and its limit $\overline Q$ satisfy the conditions in Lemma \ref{quad}.
\end{Lemma}

\subsection{Isometries}\label{app:isometries}
In this last section, we express a fundamental consequence of the previous estimate on uniform-in-$\varepsilon$ equivalence of norms. Indeed, for all exponents $a\neq 0$, the nature of the weighted Sobolev spaces changes drastically when switching between $\varepsilon>0$ and $\varepsilon=0$. For this reason, we need to embed them isometrically in the common space $\tilde H^1$ uniformly as $\varepsilon\to0$. To this aim, we can take advantage of some fundamental isometries between weighted spaces to $\tilde H^1$, which allow, by reabsorbing the weight, to obtain uniform estimates in a common space to any element in the approximating sequence.
Fixed $a\in(-\infty,1)$ and $\varepsilon\geq0$, then the map $$T^a_\varepsilon:\tilde H^1(B_1^+,\rho_\varepsilon^a(y)\mathrm{d}z)\to\tilde H^1(B_1^+)\;: u\mapsto  v=T^a_\varepsilon(u)=(\rho_\varepsilon^a)^{1/2}u $$  

is an isometry when we endow the space $\tilde H^1(B_1^+)$ with the squared norm $Q_{\rho_\varepsilon^a}$. Indeed we have:
$$\int_{B_1^+}\rho_\varepsilon^a|\nabla u|^2=Q_{\rho_\varepsilon^a}(v).$$
Is is worthwhile noticing that the family of quadratic forms $Q_{\rho_\varepsilon^a}$ is uniformly bounded (above and below) with respect to $\varepsilon\in[0,1]$.

Eventually, we remark that, similarily,  fixed $a\in(-\infty,1)$ and $\varepsilon\geq0$, then the map
\begin{equation}\label{isomsd}
\overline T^a_\varepsilon: H^1(B_1^+,\omega_\varepsilon^a(y)\mathrm{d}z)\to\tilde H^1(B_1^+)\;: u\mapsto  v=\overline T^a_\varepsilon(u)=(\omega_\varepsilon^a)^{1/2}u
\end{equation}
is also an isometry when the latter space is endowed with the squared norm $\overline Q_{\omega_\varepsilon^a}(v)$ as we have
$$\int_{B_1^+}\omega_\varepsilon^a\left(|\nabla u|^2+u^2\right)=\overline Q_{\omega_\varepsilon^a}(v).$$
Again, $\overline Q_{\omega_\varepsilon^a}$ is uniformly bounded (above and below) with respect to $\varepsilon\in[0,1]$. Once again, we remark that for these super degenerate weights Poincar\'e type inequalities do not hold true (see \cite{SirTerVit1}) and hence we can not consider only the weighted $L^2$-norm of the gradient in the equation above.

\begin{thebibliography}{99}



\bibitem{AllLSha}
M. Allen,  and H. Shahgholian, \textsl{A new boundary Harnack principle (equations with right hand side)} Arch. Ration. Mech. Anal. 234 (2019), no. 3, 1413-1444.


\bibitem{BogDyd}
K. Bogdan, B. Dyda. \textsl{The best constant in a fractional Hardy inequality}. Mathematische Nachrichten 284, 629-638 (2011).





\bibitem{CafFabMorSal} L. Caffarelli, E. Fabes, S. Mortola and S. Salsa,
\textsl{Boundary behavior of nonnegative solutions of elliptic operators in divergence form}, Indiana Univ. Math. J. 30 (1981), no. 4, 621-640.

\bibitem{CafSil1} L. Caffarelli, L. Silvestre. \textsl{An extension problem related to the fractional Laplacian}. Communications in partial differential equations 32 (2007), no. 8, 1245-1260.


\bibitem{CafSti} L. Caffarelli, P. Stinga. \textsl{Fractional elliptic equations, Caccioppoli estimates and regularity}. Ann. Inst. Henri Poincar\'e (C) Non Linear Analysis, 33, 3, (2016) 767-807.

\bibitem{DesSav}
D. De Silva, O. Savin. \textsl{A note on higher regularity boundary Harnack inequality}. DCDS-A, 35(12), (2015) 6155-6163.


\bibitem{FabKenSer}
E. Fabes, C. Kenig, and R. Serapioni. \textsl{The local regularity of solutions of degenerate elliptic equations}. Comm. Partial Differential Equations, 7(1):77-116, 1982.

\bibitem{FabJerKen1}
E. Fabes, D. Jerison, and C. Kenig. \textsl{The Wiener test for degenerate elliptic equations}. Ann. Inst. Fourier (Grenoble), 32(3):vi, 151-182, 1982.

\bibitem{FabJerKen2}
E. Fabes, D. Jerison, and C. Kenig. \textsl{Boundary behavior of solutions to degenerate elliptic equations}. In Conference on harmonic analysis in honor of Antoni Zygmund, Vol. I, II (Chicago, III., 1981), Wadsworth Math. Ser., pages 577-589. Wadsworth, Belmont, CA, 1983.

\bibitem{GilTru}
D. Gilbarg and N. S. Trudinger, \textsl{Elliptic partial differential equations of second order}, Classics in Mathematics. Springer-Verlag, 2nd edition, Berlin (1983).

\bibitem{Haj}
P. Hajlasz. \textsl{Sobolev spaces on an arbitrary metric space}. Potential Anal. 5 (1996), 403-415.





\bibitem{JavNeu}
Y. Jhaveri and R. Neumayer. \textsl{Higher regularity of the free boundary in the obstacle problem for the fractional Laplacian}. Adv. in Math. 311 (2017) 748-795.

\bibitem{JerKen}
D. Jerison and C. Kenig. \textsl{Boundary behavior of harmonic functions in nontangentially accessible domains}. Adv. in Math. 46, 80-147 (1982).




\bibitem{NorTavTerVer}
B. Noris, H. Tavares, S. Terracini, and G. Verzini. \textsl{Uniform H\"older bounds for nonlinear Schr\"odinger systems with strong competition}. Comm. Pure Appl. Math. 63 (2010) 267-302.

\bibitem{PacWei}
F. Pacard and J. Wei. \textsl{Stable solutions of the Allen-Cahn equation in dimension 8 and minimal cones}. JFA 264, 5 (2013), 1131-1167.

\bibitem{ShaYer}
H. Shahgholian and K. Yeressian, The obstacle problem with singular coefficients near Dirichlet data, Ann. Inst. H. Poincar\'e Anal. Non Lin\'eaire 34 (2017), no. 2, 293-334.


\bibitem{SirTerVit1}
Y. Sire, S. Terracini, S. Vita. \textsl{Liouville type theorems and regularity of solutions to degenerate or singular problems part I: even solutions}. arXiv:1904.02143 (2019), submitted.









\end{thebibliography}
\end{document}